\input amstex

\input xypic

\documentstyle{amsppt}
\document
\magnification=1200
\NoBlackBoxes
\nologo
\pageheight{18cm}


\bigskip

\centerline{\bf CONTINUED FRACTIONS, MODULAR SYMBOLS,} 

\smallskip

\centerline{\bf AND NON--COMMUTATIVE GEOMETRY}

\medskip

\centerline{\bf Yuri I. Manin, Matilde Marcolli}

\medskip

\centerline{Max--Planck--Institut f\"ur Mathematik, Bonn, Germany}

\bigskip

{\bf Abstract.} Using techniques
introduced by D.~Mayer, we prove an extension of the 
classical Gauss--Kuzmin theorem about the distribution 
of continued fractions,
which in particular allows one  to take into account some
congruence properties of successive convergents.
This result has an application to the Mixmaster Universe
model in general relativity.
We then study some averages involving modular symbols
and show that Dirichlet series related to modular
forms of weight 2 can be obtained by integrating certain 
functions on real axis defined in terms of continued
fractions. We argue that the quotient
$PGL(2,\bold{Z})\setminus\bold{P}^1(\bold{R})$ 
should be considered as non--commutative modular curve,
and show that the  modular complex can be seen as
a sequence of $K_0$--groups of the related
crossed--product $C^*$--algebras.


\bigskip

\centerline{\bf \S 0. Introduction and summary}

\medskip

In this paper we study the interrelation between several topics: a
generalization of the classical Gauss problem on the 
distribution of continued fractions, certain averages of modular
symbols, the properties of geodesics on modular curves, the Mixmaster
Universe model in general relativity, and the non--commutative
geometry of the quotient $PGL(2,\bold{Z})\setminus\bold{P}^1(\bold{R})$.

Our main motivation is a picture of a tower of
``non--commutative modular curves'', parameterizing two--dimensional
non--commutative tori. This ties in with the project of 
studying Stark's conjectures for real quadratic extensions of
$\bold{Q}$ via a theory of real multiplication on non--commutative
elliptic curves [Man6]. 

The traditional
algebro--geometric compactification of a modular curve $X_{G_0}= G_0
\backslash H$, for $G_0 \subset PSL(2,\bold{Z})$ a finite index
subgroup and $H$ the upper half plane, is given by the set of cusps 
$G_0 \backslash \bold{P}^1(\bold{Q})$. Our main philosophy is that
this should be replaced by the quotient $G_0 \backslash
\bold{P}^1(\bold{R})$ considered as a non--commutative space. This
point of view fits into the context of recent work [CoDS], [Man4],
[Soi].    

We support this philosophy by results of two types. First, we show
that certain exact sequences of 
[Mer], related to the modular complex introduced in [Man1], which
gives a combinatorial definition of the homology of modular curves,
can be identified with the Pimsner exact sequence for $K$--theory of our
non--commutative modular curves. Second, we demonstrate that cusps
forms of weight two for congruence subgroups (or rather their Mellin
transforms) can be obtained by integrating along the real axis certain
``automorphic series'' defined in terms of continued fractions and
modular symbols (cf.~identity (0.16)).

In a different but related perspective, we also show that the
classical definition of modular symbols can be generalized to
``limiting modular symbols'' which take into account geodesics on the
upper half plane which end at irrational points. We show that
quadratic irrationalities give rise to limiting cycles while, for
generic irrational points, there is a vanishing result in suitable
averaged sense. This result depends on the properties of a
generalization of the Ruelle transfer operator (or Gauss--Kuzmin
operator) for the shift on the continued fraction expansion, and its
spectral properties. In the case of the group $\Gamma_0(2)$ these
properties have applications to the Mixmaster Universe model in
general relativity.

\medskip

{\bf 0.1. Continued fractions.} We start by fixing the notation
which will be used throughout the paper.
Considering first $k_1,\dots ,k_n$ as independent variables, put
for $n\ge 1$
$$
[k_1,\dots ,k_n]:=
\frac{1}{k_1+\frac{1}{k_2+\dots \frac{1}{k_n}}}=
\frac{P_n(k_1,\dots ,k_n)}{Q_n(k_1,\dots ,k_n)}
\eqno(0.1)
$$
where $P_n,Q_n$ are polynomials with integral coefficients
which can be calculated inductively
from the relations
$$
Q_{n+1}(k_1,\dots ,k_n,k_{n+1})=k_{n+1}Q_n(k_1,\dots ,k_n)
+Q_{n-1}(k_1,\dots ,k_{n-1}),
$$
$$
P_n(k_1,\dots ,k_n)=Q_{n-1}(k_2,\dots ,k_n)
\eqno(0.2)
$$
It is convenient also to put formally $Q_{-1}=0, Q_{0}=1$
which is compatible with (0.2), (0.1).

\smallskip

From (0.2) one readily sees that
$$
[k_1,\dots ,k_{n-1},k_n+x_n]=
$$
$$
\frac{P_{n-1}(k_1,\dots ,k_{n-1})\,x_n+P_{n}(k_1,\dots
,k_{n})}{Q_{n-1}(k_1,\dots ,k_{n-1})\,x_n+Q_{n}(k_1,\dots ,k_{n})}=  
\left( \matrix
P_{n-1} & P_{n}\\
Q_{n-1} & Q_{n}
\endmatrix \right)\,(x_n)
\eqno(0.3)
$$
where we use the standard matrix notation for fractional 
linear transformations defining the action of $GL(2)$
on $\bold{P}^1$:
$$
z\mapsto \frac{az+b}{cz+d}=
\left( \matrix
a & b\\
c &  d
\endmatrix \right)\,(z)
$$

\smallskip

If $\alpha \in (0,1)$ is an irrational number,
there is a unique sequence of integers $k_n(\alpha )\ge 1$
such that $\alpha$ is the limit of $[\,k_1(\alpha ),\dots ,k_n(\alpha )\,]$
as $n\to\infty$. Moreover, there is a unique
sequence $x_n(\alpha )\in (0,1)$ such that
$$
\alpha = [\,k_1(\alpha ),\dots ,k_{n-1}(\alpha ), k_n(\alpha )+x_n(\alpha )\,]
$$
for each $n\ge 1$. Rational numbers in $(0,1]$ can be accommodated by allowing
finite sequences of $k_i\ge 1$ complemented by zeroes, and similarly
for $x_n(\alpha )$.

\smallskip

We can specialize (0.3) at the point $\alpha$ and get by induction
$$
\alpha =
\left( \matrix
0 & 1\\
1 & k_1(\alpha )
\endmatrix \right)\ \dots \
\left( \matrix
0 & 1\\
1 & k_n(\alpha )
\endmatrix \right)\,(x_n(\alpha ) ).
\eqno(0.4)
$$
Put 
$$
p_n(\alpha ):=P_n(k_1(\alpha ),\dots ,k_n(\alpha )),\
q_n(\alpha ):=Q_n(k_1(\alpha ),\dots ,k_n(\alpha ))
$$
so that $p_n(\alpha )/q_n(\alpha )$ is the sequence
of convergents to $\alpha$. Denote also 
$$
g_n(\alpha ):=\left( \matrix
p_{n-1}(\alpha ) & p_n(\alpha )\\
q_{n-1}(\alpha ) & q_n(\alpha )
\endmatrix \right)
$$
\smallskip

For further use, we reproduce a description of the total set
of matrices  $g_n (\alpha )$. Put
$$
\roman{Red}_n=\left\{ \,
\left( \matrix
0 & 1\\
1 & k_1
\endmatrix \right)\ \dots \
\left( \matrix
0 & 1\\
1 & k_n
\endmatrix \right)\,|\,k_1,\dots , k_n\ge 1;\, k_i\in\bold{Z}\right\},
\eqno(0.5)
$$
and $\roman{Red}:=\cup_{n\ge 1}\,\roman{Red}_n \subset GL(2,\bold{Z}).$

\smallskip

In the terminology of [LewZa1], $\roman{Red}_n$ consists
of reduced matrices of length $n$. The following properties
of $\roman{Red}$ are proved in that paper:

\smallskip

(i) {\it A matrix in $GL(2,\bold{Z})$ is reduced iff it has
non--negative entries which are non--decreasing downwards
and to the right.}

\smallskip

(ii) {\it The length $l(g)$ of a reduced matrix $g$ and its representation
in the form (0.5) are uniquely defined.}

\smallskip

(iii) {\it Reduced matrices are hyperbolic and have two distinct
fixed points on the real axis. 
Every conjugacy class $g$ of hyperbolic matrices  in
$GL(2,\bold{Z})$ contains reduced representatives.
They all have the same length $l(g)$, and there are
exactly $l(g)/k(g)$ of them where $k(g)$ is the maximal
integer such that $g=h^{k(g)}$ for some $h$.}

\medskip

{\bf 0.1.1. Generalized Gauss problem.} Consider
a subgroup of finite index $G\subset GL(2,\bold{Z})$
and the coset space $\bold{P}= GL(2,\bold{Z})/G$ with the left transitive
action of  $GL(2,\bold{Z})$ on it: 
$$
t\mapsto \left( \matrix
a & b\\
c &  d
\endmatrix \right)\,(t)
$$
For any $x\in [0,1],\,t\in\bold{P},n\ge 0$ put
$$
m_n(x,t):= \roman{measure\ of\ the\ set}\
\{\alpha\in (0,1)\,|\,x_n(\alpha )\le x,\,g_n(\alpha)^{-1}(t_0)=t\} 
\eqno(0.6)
$$
where $t_0$ is the base point of $\bold{P}$, the coset of $G$. 
Notice that $x_n(\alpha )=g_n(\alpha)^{-1}(\alpha )$
so that $m_n$ is essentially the pullback of the Lebesgue
measure on $(0,1)\times \bold{P}$ with respect to the
operator $g_n(\alpha )$ acting upon $\alpha$ and $t$ simultaneously.
Notice also that $g_n(\alpha )^{-1}$ is another notation for the
$n$--th power of the shift operator 
$$ 
T:\, (\alpha ,t)\mapsto \left( \frac{1}{\alpha}-\left[
\frac{1}{\alpha} \right], \left( 
\matrix -[1/\alpha ] &  1\\ 1 & 0 \endmatrix \right)\,(t)
\right) 
$$

\smallskip

The first result of this paper is the following
generalization of the Gauss--Kuzmin--L\'evy formula:

\smallskip

\proclaim{\quad 0.1.2. Theorem} Assume that $\roman{Red}\,(t)=
\bold{P}$ for each $t\in \bold{P}$ (transitivity condition), with
$\roman{Red}\,(t):= \{ gt | g\in \roman{Red} \subset GL(2,\bold{Z})
\}$. Then the limit
$m(x,t) = \roman{lim}_{n\to\infty}\,m_n(x,t)$ exists
and equals
$$
m(x,t)=\frac{1}{|\bold{P}|\,\roman{log}\,2}\,\roman{log}\,(1+x).
\eqno(0.7)
$$
\endproclaim
\smallskip

A proof is given in \S 1 below. It starts with a straightforward 
generalization of the Gauss--Kuzmin
inductive expression of $m_{n+1}$ through $m_n$. To write it down,
consider first the following sets: for $y\in (0,1), s\in \bold{P}$ put
$$
M_n(y,s):= \{\beta \in (0,1)\,|\,x_n(\beta )\le y,\, g_n(\beta )^{-1}t_0=s\}.
$$
Then we have, using (0.4) and neglecting the rationals which have measure zero:
$$
M_{n+1}(x,t)=\coprod_{k=1}^{\infty}
\left\{M_n \left(\frac{1}{k},
\left( \matrix
0 & 1\\
1 & k
\endmatrix \right)\,(t)\right) - 
M_n \left(\frac{1}{x+k},
\left( \matrix
0 & 1\\
1 & k
\endmatrix \right)\,(t)\right)\right\} .
$$
Therefore,
$$
m_{n+1}(x,t)=\sum_{k=1}^{\infty}
\left\{m_n \left(\frac{1}{k},
\left( \matrix
0 & 1\\
1 & k
\endmatrix \right)\,(t)\right) - 
m_n \left(\frac{1}{x+k},
\left( \matrix
0 & 1\\
1 & k
\endmatrix \right)\,(t)\right)\right\} .
\eqno(0.8)
$$
Derivating (0.8) in $x$, we get the following equation
for the densities, which introduces an important operator $L$,
the generalized Gauss--Kuzmin operator:
$$
m_{n+1}^{\prime}(x,t)= (Lm_n^{\prime})(x,t):=
\sum_{k=1}^{\infty}
\frac{1}{(x+k)^2}
m_n^{\prime} \left(\frac{1}{x+k},
\left( \matrix
0 & 1\\
1 & k
\endmatrix \right)\,(t)\right) .
\eqno(0.9)
$$

In the classical case, when $|\bold{P}|=1$, at least
four ways to deduce (0.7) from (0.9) are known:
two ``elementary'' deductions going back to R.~Kuzmin
and P.~L\'evy respectively (see e.g. [Sch]), and two functional
analytic methods, based upon the spectral analysis
of the operator $L$ formally defined by (0.9):
see [Ba], [BaYu] and [May1].

\smallskip

Clearly, the limiting measure must be an eigenfunction
of $L$ corresponding to the eigenvalue $1$,
and both analytic proofs show the convergence
of $L^nm_0^{\prime}$ to this limiting measure
by establishing that $1$ is the multiplicity one
eigenvalue with maximal modulus for the extension
of $L$ to an appropriate function space.
K.~Babenko realizes $L$ as a self--adjoint operator
in a Hilbert space, whereas D.~Mayer works
in the context of nuclear and trace class operators
in Banach spaces.

\smallskip

Of these four proofs, we were able to generalize
to our context only Mayer's method. Babenko's
representation seems to be inadequate for proving convergence.
However, it might still be useful for numerical calculations,
and we present it in \S 1.3.

\smallskip

The operator $L$ and its deformation $L_s$ ((1.1) below)
were introduced and studied also in a recent
paper [ChMay], of which we became aware only after the
first draft of this paper was written.

\medskip

{\bf 0.2. Modular curves and geodesics.} Our study of
the generalized Gauss measure (0.7) was motivated
by the relationship between continued fractions
and one--dimensional homology of modular curves.
We will start by recalling the basic features of
this relationship in the form given in [Man1]; see
also [Mer] for additional information.

\smallskip

Let $G_0$ be a subgroup of finite index in
the fractional linear group
$PSL(2,\bold{Z})=SL(2,\bold{Z})/(\pm 1).$
It determines the noncompact modular
curve $G_0\setminus H$ (where $H$ is the upper complex half--plane).
This curve admits a smooth compactification by a finite number of
cusps which are in a natural bijection 
with the set $G_0\setminus\bold{P}^1(\bold{Q}).$
Let $X_{G_0}$ (or more precisely, $X_{G_0}(\bold{C})$)
denote this compactification, $\varphi$ the respective
covering map.

\smallskip

For any two points $\alpha ,\beta\in \overline{H}:=H\cup
\bold{P}^1(\bold{Q})$ we can define a real homology
class (``modular symbol'') $\{\alpha ,\beta\}_{G_0}\in
H^1(X_{G_0},\bold{R})$ by integrating the lifts of differentials
$\omega$ of the first kind on $X_{G_0}$ along the geodesic path
connecting $\alpha$ to $\beta$:
$$
\int_{\{\alpha ,\beta\}}\omega :=\int_{\alpha}^{\beta} \varphi^*(\omega ).
\eqno(0.10)
$$
Modular symbols satisfy
the basic additivity and invariance properties:
$$
\{\alpha ,\beta\}_{G_0}+\{\beta ,\gamma\}_{G_0}=\{\alpha ,\gamma\}_{G_0},
$$ 
$$
\forall g\in G_0,\ \{g\alpha ,g\beta\}_{G_0}=\{\alpha ,\beta\}_{G_0}.
\eqno(0.11)
$$

The integrals (0.10) can be related to  finite, stably periodic, or
general infinite continued fractions, depending on the arithmetical
nature of the ends $\alpha ,\beta$. We will briefly treat
these three cases separately. 

(i) {\it Finite continued fractions.} Assume first that
$\alpha ,\beta$ are cusps. It is known that in this case
the modular symbol represents
a rational homology class (Manin--Drinfeld's theorem).

\smallskip

By additivity, it suffices to look at the integrals of the form
$\int_0^{\alpha},\,\alpha\in \bold{Q}$. Let 
$$
g_k:=\left( \matrix p_{k-1}(\alpha )& p_{k}(\alpha )\\
q_{k-1}(\alpha )& q_{k}(\alpha ) \endmatrix \right),\
k=1,\dots ,n,\ \alpha =\frac{p_n(\alpha)}{q_n(\alpha )}.
$$
Then, again by additivity, we have
$$
\int^{\alpha}_0 \varphi^*(\omega )=
\sum_{k=1}^n\, \int_{p_{k-1}(\alpha )/q_{k-1}(\alpha )}^{p_k(\alpha )/q_k(\alpha )} \varphi^*(\omega )=
-\sum_{k=1}^n\, \int_{\{g_k (0),g_k(i\infty)\}}\omega .
\eqno(0.12)
$$
Finally, in view of (0.11), the $k$--th integral in (0.12)
depends only on the class of $g_k$ in $G_0\setminus PSL (2,\bold{Z}).$

\smallskip

Thus, (0.12) establishes a connection between the distribution
of modular symbols $\{\alpha ,\beta\}_{G_0}$ and the distribution
of pairs of consecutive convergents to $\alpha$ and
$\beta$ in $\bold{P}_0:=G_0\setminus PSL (2,\bold{Z}).$ 

\smallskip

Notice
that there is a slight discrepancy
with 0.1.1 where we dealt instead with $G\subset GL (2,\bold{Z})$
and $GL (2,\bold{Z})/G.$ In order to reduce the current problem
to the former one, we can first replace $G_0$ by its lift
to $SL (2,\bold{Z})$ and then denote by $G$ the subgroup
generated by this lift and $\left( \matrix -1 & 0\\ 0 & 1
\endmatrix \right)$. We will have a natural identification
$\bold{P}_0= G\setminus GL(2,\bold{Z})$, and this set
in turn identifies with $\bold{P}$ in 0.1.1 under
the map $g\mapsto g^{-1}$ which is implicit in (0.6).

\smallskip

If $\alpha$ and/or $\beta$ in (0.10) is real irrational,
the integral diverges at this end.
In this case,
it is natural to study its asymptotic behavior.
We will define the ``limiting modular symbol''
by the following
expression whenever it makes sense:
$$
\{\{* ,\beta\}\}_{G_0}:=\roman{lim}\,
\frac{1}{T(x,y)}\,\{x,y\}_{G_0}\in H^1(X_{G_0},\bold{R}).
\eqno(0.13)
$$
Here $x,y\in H$ are two points on the geodesic joining
$\alpha$ to $\beta$, $x$ is arbitrary but fixed,
$T(x,y)$ is the geodesic distance between them, and
the limit is taken as $y$ tends to $\beta$.
In \S 2 we will prove that if the limit does exist, 
it depends neither on $x$ nor
on $\alpha$, hence the first argument is replaced by $*$,
whereas the double curly brackets remind about the limit.
We will discuss the relation  of this symbol to continued
fractions
in two situations.

\smallskip

(ii) {\it Stably periodic continued fractions.} 
Among geodesics with two irrational ends
there is an important subclass consisting
of geodesics that connect two fixed points of
a hyperbolic element in $G_0$. More precisely,
any hyperbolic $g\in G_0$ has two fixed points
$\alpha^{\pm}$, repelling and attracting,
on the real line. Let $\Lambda_g^{\pm}$
be the respective eigenvalues, $0<\Lambda_g^{+}<1.$
The oriented geodesic in $H$ connecting $\alpha^{-}_g$
to $\alpha^{+}_g$ is $g$--invariant, and the action
of $g$ induces on it the shift by the geodesic
distance $\lambda (g):=\roman{log}\,\Lambda^-_g$.

\smallskip

For any point $x$ on this geodesic, the image of its
segment $[x,gx]$ is a parameterized closed loop
on $X_{G_0}$ missing the cusps. (The supporting
set--theoretic loop is ran over exactly $k(g)$
times where $k(g)$ is the maximal $k$ such
that $g$ is a $k$--th power in $G_0$).
The homology class of this loop is $\{0,g(0)\}$.
When we integrate from a fixed $x$ to $g^nx$, we 
run over the parameterized loop
$n$ times, the geodesic length
of the path is $n\,\lambda (g)$, and its
homology class is $n\,\{0,g(0)\}$.
Therefore the limiting modular symbol (0.13) exists and equals
$$
\{\{*,\alpha^+_g\}\}=\frac{\{0,g(0)\}}{\lambda (g)}.
\eqno(0.14)
$$
\smallskip

The most important generating function
for closed geodesics is the Selberg zeta function.
However, it encodes only 
the lengths of closed geodesics in the hyperbolic metric.
The usual modular symbol in the numerator of (0.14)
depends only on the class of $g$ modulo $[G_0,G_0]$
and is additive in $g$
(see [Man1], Prop. 1.4). Perhaps, one can construct
a generating function for (0.14) as a combination
of Selberg's zetas with abelian characters.

\smallskip

The usual Selberg's zetas were studied
in [May1], and then in [LewZa1] for $GL$ and $SL$ separately.
It turned out that they could be represented
as Fredholm determinants $\roman{det}\,(1-L_s)$
and $\roman{det}\,(1-L_s^2)$ respectively.
Literally the same is true in our generalized
setting, when subgroups $G$ or $G_0$ are introduced.
This is proved in [ChMay], and we supply a brief
discussion of this in \S 3.

\smallskip

As we have briefly explained, the distribution of 
continued fractions and modular symbols at cusps is encoded
in the eigenvalues of $L_1$, $1$ being the dominant
value producing the distribution (0.7).

\smallskip

From the identity $\roman{det}\,(1-L_s)=Z(s)$
it follows that the zeroes of $Z(s)$ are exactly those
values for which the deformed operator $L_s$ 
has  eigenvalue $1$.

\smallskip

J.~Lewis and D.~Zagier produced also an in--depth study of
the respective eigenfunctions for
the full modular groups $GL(2,\bold{Z})$ 
and $SL(2,\bold{Z}).$

\smallskip

(iii) {\it General infinite continued fractions.}
For this case, we will prove in \S 2 two results.
Namely, we will establish that with an additional
assumption the limiting modular symbol
(0.13) exists only in a weak sense and and is zero.

\smallskip

\proclaim{\quad 0.2.1. Theorem} 
Assume that $\roman{Red}(t)=\bold{P}_0$
for each $t\in\bold{P}_0$. 
Then (0.13) weakly converges to
zero.
\endproclaim

\smallskip

For the precise description of the sense in which
this convergence holds, see \S 2.3 and in particular,
(2.21). Our prof is based upon Theorem 0.1.2.

\smallskip

This vanishing can be compared with a well known interpretation
of the Selberg trace formula for compact surfaces:
quantum mechanical averages for
the geodesic flow can be calculated as if this
flow were classically concentrated on closed geodesics.

\smallskip

The transitivity assumption for $\roman{Red}$
will be checked in \S 2 in the case $G_0=\Gamma_0(N).$

\smallskip

We show that the case $N=2$ has a nice little application
to the study of the $t=0$ singularity of the Bianchi IX
model in general relativity.

\smallskip

The next result of \S 2 concerns a series of averaging
formulas of a different kind.
Drawing on a lemma of P.\,L\'evy,
we will explain in \S 2 how to calculate
averages (over $[0,1]$) of some functions of $\alpha$ 
defined by the sums over all pairs of consecutive convergents of
$\alpha$. Here we will state an interesting particular
case, providing averages of weighted
modular symbols.
 
\medskip

Fix a prime number $N>0$ and put $G_0=\Gamma_0(N)$.
We will assume that the genus of $X_{G_0}=X_0(N)$ is $\ge 1$,
otherwise our identities become trivial.
Consider a $\Gamma_0(N)$--invariant differential
$\varphi^*(\omega )$ on $H$ which is a cusp eigenform for all
Hecke operators and denote by $L_{\omega}^{(N)}(s)$
(resp. $\zeta^{(N)}(s)$)
its Mellin transform (resp. Riemann's zeta) with omitted Euler $N$--factor.
More precisely, the coefficients of $L_{\omega}^{(N)}(s)$
are Hecke eigenvalues of $\varphi^*(\omega )/dz$.

\medskip

\proclaim{\quad 0.2.2. Theorem} 
We have for $\roman{Re}\,t>0$:
$$
\int_0^1 d\alpha \sum_{n=0}^{\infty}
\frac{q_{n+1}(\alpha )+q_n(\alpha )}{q_{n+1}(\alpha )^{1+t}}\,
\int_{\{0,q_n(\alpha )/q_{n+1}(\alpha )\}}\omega  =
$$
$$
\left[ \frac{\zeta (1+t)}{\zeta (2+t)}
-\frac{L_{\omega}^{(N)}(2+t)}{\zeta^{(N)} (2+t)^2}
\right]\,\int_0^{i\infty}\varphi^*(\omega ).
\eqno(0.16)
$$
\endproclaim

\medskip

Our calculation of averages of modular symbols like 
(0.16) and various generalizations in \S 2.1--2.2
point a way towards understanding what function theory on
non--commutative modular curves may be used in order to recover the  
theory of modular forms on the upper half plane. In fact, (0.16)
represents Mellin transforms of weight two cusp forms in terms of the
quantities that can be defined entirely 
in terms of the {\it boundary} of the moduli space,
and not the traditional integrals along
geodesics: in fact, the integral in the left hand
side of (0.16) is taken along the real axis.
 
\smallskip

For us this boundary is $\bold{P}^1(\bold{R})$, and
not $\bold{P}^1(\bold{Q})$ as in the traditional algebro--geometric
compactification, and it is exactly the consideration
of this boundary that leads us into the land of non--commutative
geometry: since modular symbols, Hecke operators and
continuous fractions all can be expressed in terms of the
(noncommutative) geometry of the boundary, so our $L$--series as well can.
(The appearance of $\omega$ can be avoided since
we could consider only eigenvalues of Hecke
operators acting on the modular complex.)

\medskip

{\bf 0.3. Relations with non--commutative geometry.} As is well--known,
the quotient $PGL(2,\bold{Z})\setminus \bold{P}^1(\bold{R})$
can be identified with the space of classes of continued fractions
modulo the equivalence relation ``{\it $k_{n+n_0}(\alpha
)=k_{n+n_1}(\beta )$ for some $n_0,n_1$ and all $n$}''. Classical
results on various averages like 
$$
\roman{lim}_{n\to\infty}\,n^{-1}\sum_{i=1}^n f(k_i(\alpha ))
$$
state that these averages are almost everywhere constant
functions on this space.

\smallskip

On the other hand, this space and its finite 
coverings corresponding to subgroups 
$G\subset PGL(2,\bold{Z})$ constitute the boundary
of the analytic moduli stack of elliptic curves whose irrational
points are invisible in algebraic geometry (only cusps
admit classical algebraic interpretation). According to 
the emerging general philosophy, this boundary is a bridge
to the world of non--commutative geometry.
In particular, the geometric objects parameterized
by this boundary, which are two--dimensional non--commutative tori
modulo Morita equivalence,
can be treated as limiting elliptic curves.
For some explanations, see [CoDS], [Soi],
[Man4].

\smallskip

Accordingly, the boundary itself
should be considered as (a tower of) ``non--commutative 
modular curves'' in Connes' spirit.
The modular complex introduced
in [Man1] and further studied in [Mer], [Gon] and
other papers,  provides a combinatorial definition
of the homology of the modular tower. In \S 4
we show that essentially the same complex calculates
$K$--theory of the crossed--products describing
the non--commutative boundary modular tower.

\smallskip

This viewpoint presents in the new light also 
the identity (0.16) and its generalizations
considered in \S 2.1--2.2. Namely, it demonstrates that
at least a part of the theory
of modular forms in the upper half
plane can be recast 
as the study of averages of certain functions defined on the
boundary $\bold{R}$
as sums of the type (2.1). Their behavior with respect
to fractional linear transformations is not modular
in the traditional sense, but their expression
via pairs of successive denominators can be seen as remnants of modularity. 

\smallskip

For another family of similar phenomena, see [Za] and [LawZa]. 

\smallskip

It is interesting to remark in this context that
the Gauss--Kuzmin operator $L_s$ vaguely looks like
a ``Hecke operator at the arithmetical infinity'',
and has some properties that might
be expected of such an operator.

\smallskip

As a final remark, in [Man2]
it was shown that after a choice of Schottky
uniformization, the Arakelov geometry of a complex
curve $X$ at arithmetical infinity can be described
in terms of the hyperbolic geometry of geodesics
not on $X$ itself, but rather in the hyperbolic
handlebody having $X$ as the boundary at infinity. It would be 
interesting to clarify the statistical
aspects of the closed and infinite geodesics in the same
vein as above and to relate them to Arakelov geometry.
We hope to return to this question later.

\bigskip

{\bf Acknowledgment.} We thank Dieter Mayer, Victor Nistor, and Don
Zagier for useful conversations. The second author is partially
supported by Sofja Kovalevskaya Award.

\newpage

\centerline{\bf \S 1. Gauss--Kuzmin operator}

\medskip

{\bf 1.1. Operator $L_s$.} Consider the operator $L_s$
acting on functions of two variables $(x,t)$ and formally
given by 
$$
(L_sf)(x,t):=\sum_{k=1}^{\infty}
\frac{1}{(x+k)^{2s}}
f \left(\frac{1}{x+k},
\left( \matrix
0 & 1\\
1 & k
\endmatrix \right)\,(t)\right) .
\eqno(1.1)
$$
The variable $x$ here varies in a subset of
$\bold{C}$ stable with respect to all maps
$x\mapsto (x+k)^{-1},\, k=1,2,\dots$. In our
context this subset will always contain $[0,1]$. The variable
$t$ belongs to a finite $GL(2,\bold{Z})$--set $\bold{P}$
endowed with a base point $t_0$.
The parameter $s$ here is real and $>1/2$. We are mostly
interested in the case $s=1$. In \S 3 we will allow
$s$ to take complex values. Finally, $f$ will vary
in a linear space of functions stable with respect
to $L_s$ and containing the function
$m_0^{\prime}(x,t)=\delta_{t,t_0}$ (cf. 0.1.1).

\smallskip

For the proof of Theorem 0.1.2, we want to create a functional
analytic context in which the machinery of Krasnoselskii's
theorem as stated in [May1], 7.25, becomes applicable.
To this end we make the following choices,
slightly generalizing Mayer. (They were also made
in [ChMay]).

\smallskip

(i) {\it Definition domain.} It will be $\bold{D}\times \bold{P}$
where $\bold{D}:= \{z\in\bold{C}\,|\,|z-1|<3/2\}.$
We will call the subsets $\bold{D}\times \{t\}$ {\it sheets}.
Notice that each map $z\mapsto (z+k)^{-1}$ transforms
$\bold{D}$ strictly into itself.

\smallskip

(ii) {\it Functional spaces}. We shall consider the
complex Banach space $B_{\bold{C}}:=V_{\bold{C}}(\bold{D}\times \bold{P})$ consisting of functions
holomorphic on each sheet and continuous on its boundary. We shall
also consider the real Banach space $B:=V(\bold{D}\times
\bold{P})$ of functions holomorphic on each sheet and continuous
on its boundary, real at the real points of each sheet. Both spaces
are endowed with the supremum norm.
These spaces obviously contain $\delta_{t,t_0}$
and are stable with respect to $L_s$ for real $s>1/2.$
The space $B_{\bold{C}}$ is also stable with respect
to $L_s$ with $\roman{Re}\,s>1/2.$

\smallskip

(iii) {\it Positive cone.} Denote by $K\subset B$ the cone
consisting of functions taking non--negative values
at real points of each sheet. We have $K\cap -K={0}$
($K$ is proper)
because a nonzero analytic function cannot vanish on an interval.
We also have $B=K-K$ ($K$ is reproducing) because $f= (f+r)-r$,
and if $r$ is large and positive, $f+r, r\in K.$
Finally, functions positive at all real points of all sheets
form the interior of ${K}$.

\smallskip

We write $f\le g$ if $f-g\in K.$

\smallskip

(iv) {\it $L_s$ is $K$--positive.} This means that
$L_s(K)\subset K$ which is obvious. 

\medskip

\proclaim{\quad 1.1.1. Lemma} Assume that $\bold{P}$
contains no proper invariant subsets with respect
to the operators $\roman{Red}$
(see (0.5)).

\smallskip
Then for each nonzero $f\in K$
there exist two real positive constants $a,b$ and
an integer $p\ge 1$ such that $a\le L_s^pf\le b.$

\endproclaim

\smallskip

{\bf Proof.} The upper bound is trivial. Assume 
that for some $f$ and all $p$ the lower bound is zero. This
means that for each $p\ge 1$, $L_s^pf$ vanishes
at some point $(x_p,t_p)$ with $x_p$ real in the closure of $\bold{D}.$
Since in (1.1) all summands are non--negative at real points,
when $f\in K$,
we see that  $f$ must vanish at all points contained
in the set $\cup_{p\ge 1}\roman{Red}_p(x_p,t_p)$.

\smallskip

From our assumption it follows that for  some $q$ and any $t\in \bold{P}$,
$\roman{Red}_q (t)=\bold{P}$. By downward induction we deduce first
that for any $p$, $L^p_sf$ has real zeroes on all sheets. Then, again
by downward induction, 
one sees that for each $t$  there exists  a sequence
of integers $q_n\to\infty$ and real points $y_n$
in $\bold{D}$ such that $f(x,t)$ vanishes
at all $x\in \cup_n\roman{Red}_{q_n}(y_n)$. But the intersection
of the latter set with $[0,1]$ is dense in $[0,1]$. A nonvanishing
holomorphic function cannot have as many zeroes.
This contradiction proves our assertion.

\medskip

\proclaim{\quad 1.1.2. Lemma} $L_s:\,B_{\bold{C}}\to B_{\bold{C}}$ 
is a nuclear operator
of order zero, in particular compact and trace class for
$\roman{Re}\,s>1/2.$
\endproclaim  

\smallskip

{\bf Proof.} The reasoning is the same as in [May1], and we
only sketch it. Denote the $k$--th summand in (1.1)
by $(\pi_{s,k}f) (x,t)$.
Each $\pi_{s,k}$ is nuclear, and $\sum_k \|\pi_{s,k}\|$
converges for $\roman{Re}\,s>1/2$. In fact, the spectrum
of  $\pi_{s,k}$ can be easily calculated.
Let $z_k$ be the unique fixed point of
$\gamma_k:=\left( \matrix
0 & 1\\
1 & k
\endmatrix \right)$
in $\bold{D}$, and $\mu_i^{(k)}$ the spectrum of the permutation
induced by this matrix on $\bold{P}$. Then
the spectrum of $\pi_{s,k}$ is $\{(-1)^n(z_k+k)^{-2(s+n)}\mu_i^{(k)}\}$,
$n\ge 0.$

\medskip

{\bf 1.2. Proof of the generalized Gauss--Kuzmin Theorem.}
We have now checked all the conditions for
the applicability of the Theorem 7.25 of [May1]; see also
[KraLS] for more details.

\smallskip

Using this theorem, we conclude that
there exists exactly one eigenfunction of norm one $f_s$ of $L_s$
in the interior of $K$ and its eigenvalue $\lambda_{0,s}$ is positive
and simple. All other eigenvalues have strictly lesser
modulus. For any $f\in B$ and $\epsilon >0$, we have
$$
L_s^n f= \roman{const}\,\lambda_{0,s}^n\,f_s +O(c(\epsilon)(q
+\epsilon)^n \lambda_{0,s}^n) 
\eqno(1.2)
$$
as $n\to \infty$,
where $q=q(L_s) <1$ is the spectral margin of $L_s$.

\smallskip

In the case $s=1$ we know a positive eigenfunction:
this is Gauss' density $\dfrac{1}{x+1}$, independent
of $t$. The normalization constant is straightforward.

\smallskip

This argument completes the proof of the generalized
Gauss--Kuzmin theorem. We will, however, provide some more
details about the deduction of (1.2), because this
technique can be useful also in the treatment of lower eigenvalues.

\smallskip

The basic result that ensures the existence of eigenfunctions in
certain invariant cones is the following ([KraLS], Theorem 9.2):

\smallskip

{\it If $K$ is a cone in a real Banach space $B$ satisfying
$\overline{K-K}=B$ and $L$ is a compact operator with $LK \subseteq
K$, and with positive spectral radius $r(L)$, then $r(L)$ is an
eigenvalue of $L$ with a corresponding eigenfunction in $K$.}

\smallskip

We also recall some results which enable us to establish that the top
eigenvalue of an operator with an invariant cone is simple, and when the
rest of the spectrum is of strictly smaller absolute value.

Following [May1], we say that an operator $L$ is $u$--bounded, with respect
to a function 
$u\in K$, if for any $f\in K$ there exists some $n >0$ and $a,b>0$
such that
$$ 
au \le L^n f \leq bu 
$$
(Since we allow a power of $L$, this is a weaker 
definition of $u$--boundedness than the one
on p.~110 of [KraLS], cf.~their remark on p.~111.)

\smallskip

Lemma 1.1.1 shows that $L_s$ is $u$--bounded
with respect to the constant
function $u(x,t)=1$. The lower bound guarantees the positivity of the
spectral radius $r(L_s)$ and hence the applicability 
of the Theorem 9.2 of [KraLS]. This fact
follows from Lemma 9.2 of [KraLS].

\smallskip
 
We then have the following result ([KraLS], Theorem 11.1).
Assume that the cone $K$ is reproducing and the $K$--positive operator $L$ is
$u$--bounded. Assume moreover that $L$ has an eigenvalue $\lambda_0 >0$
with an eigenvector $f\in K$. Then the eigenvalue $\lambda_0$ is
simple. A proof of this fact can be obtained using Lemma 11.1 of
[KraLS]. This is actually a simple result of
linear algebra which uses only the fact that $f$ is an interior point
of the cone (nowhere vanishing), and that iterates of the operator
$L$ map boundary points of the cone different from $\{ 0 \}$ to interior
points of the cone. 

\smallskip

Furthermore, Theorem 11.4 of [KraLS]
shows that with the same hypothesis as in the previous result,
every eigenvalue $\lambda$ of $L$ different from $\lambda_0$ satisfies
$| \lambda | < \lambda_0$.
This result follows from the observation that if the operator $L$ is
$u$--bounded then it is also $f$--bounded, where again $f$ is the
eigenfunction in $K$ with eigenvalue $\lambda_0$. Then, if $h$ is an
eigenfunction with eigenvalue $\lambda$ the estimate
$$ -\alpha (\lambda_0 -\epsilon)f \leq \lambda h \leq \alpha(\lambda_0
-\epsilon) f $$
for some $\epsilon >0$ follows easily, where $\alpha>0$ is the
smallest positive number such that $-\alpha f \leq h \leq \alpha f$
is satisfied. This gives $|\lambda| \leq \lambda_0-\epsilon$.

\smallskip

Since in our case we know that the operator $L_s$ is compact, the
previous result implies that all the other eigenvalues $\lambda$
(hence all the points in the spectrum of $L_s$) satisfy the estimate
$|\lambda | < q\lambda_0$, for some $q<1$. 

\smallskip   

Finally, we have a result on the convergence of iterates, cf. Theorem
15.4 of [KraLS].

\smallskip 

The cone $K$ in our case contains some ball
of positive radius. In this case, Theorem 9.11 on
p.~97 of [KraLS] ensures that the adjoint operator $L^*$ acting
on the dual Banach space $B^{\prime}$ has an eigenfunctional $f^*$ in the
adjoint cone $K^*$ of linear $K$--positive functionals,
with eigenvalue $\lambda \leq r(L)$ where $r(L)$ is the spectral radius
of $L$. In our case $L_s$ has an eigenvalue $\lambda_0 = r(L_s)$ and a
corresponding eigenfunction $f$ which is an interior point of the cone
$K$. Thus, if $L^* f^* = \lambda f^*$, for a non-trivial $f^*\in K^*$,
and if moreover $f^*(f) >0$, then $\lambda =\lambda_0$, because
$$ 
\lambda_0 f^*(f) = f^*(\lambda_0 f) = f^*(Lf) = (L^* f^*)(f) =
\lambda f^*(f). 
$$

\smallskip

Assume that the operator $L$ has a simple
eigenvalue equal to the spectral radius, $\lambda_0=r(L)$, and the
remaining part of the spectrum lies in the disk $| \lambda | < q\,r(L)$
for some $q<1$. Let $f\in K$ be the eigenfunction of the eigenvalue
$\lambda_0$. Let $f^*$ be an eigenfunctional for $L^*$ in
$K^*$, with eigenvalue $\lambda_0$ as above, satisfying $f^*(f)=1$.
Then the sequence of iterates
$$ f_{n+1} = L f_n $$
converges to the eigenfunction $f$ in the following sense. 

\smallskip

Define the operators $Uh :=f^*(h) f$, and $U^\perp
h := h - f^*(h) f$. 
We have
$$ 
\lim_n \frac{\| U^\perp f_n \|}{\| U f_n \|} =0 
$$
and the rate of convergence is estimated by
$$
\frac{\| U^\perp f_n \|}{\| U f_n \|} \leq c(\epsilon) (q
+\epsilon)^n  \frac{\| U^\perp f_1 \|}{\| U f_1 \|}, 
$$
for arbitrarily small $\epsilon >0$. In other words, the iterates
converge as fast as a geometric progression with ratio arbitrarily
close to the spectral margin $q=q(L)$, cf. \S 15.2 of [KraLS], and in
particular Theorem 15.3, where a more refined estimate of the
coefficient $c(\epsilon )$ is also given. 

\smallskip

Sometimes similar techniques may be applied to the study of the second
eigenvalue, by supplying a suitable real Banach space with an
invariant cone for the operator $h \mapsto L_s h -\lambda_{0,s} f^*
(h) f$, cf. [May1]. 

\medskip

As an example, we will now show that the condition of Lemma 1.1.1
holds for the congruence subgroup $\Gamma_0(N).$ For $N=2$,
the generalized Gauss--Kuzmin Theorem has a nice application to
the dynamical system arising in the general relativity,
the so called ``Mixmaster Universe''.

\medskip

\proclaim{\quad 1.2.1. Proposition} Let $G$ be the subgroup generated
by the lift of $\Gamma_0(N)$ and the sign change.
Then $\roman{Red}_3(t)=\bold{P}$
for any $t$.
\endproclaim

\smallskip

{\bf Proof.}  In fact, elements of $\bold{P}$ can be thought of
as points of the projective line over $\bold{Z}/N$,
that is, formal quotients of residues $\roman{mod}\,N$ that
can be represented by pairwise prime integers ([Man1]).
Moreover, this encoding can be chosen compatible with
the usual action of $GL(2)$ upon $\bold{P}^1$.

\smallskip

Let us break these points into three groups:

\smallskip

(I) $\{u/1\,| \,(u,N)=1\}.$

\smallskip

(II) $\{du/1\,|\,d/N, d>1, (u,N)=1\}.$ 

\smallskip

(III) $\{1/du\}\,|\,d/N, d>1, (u,N)=1$.

\smallskip

Let us say that $s$ can be obtained from $t$ in one
step, if $s\in \roman{Red}_1 (t)$, that is,
$s=(t+k)^{-1}$ for some $k.$ The following
statements are straightforward, and taken together,
prove our claim.

\smallskip

From any single element of I\,$\cup$\,II we can obtain in one step 
all elements of the set I\,$\cup$\,III.

\smallskip

From any single element of III one can obtain in one step
an element of II, by adding zero and inverting.
Hence from the total III one can obtain the total II in one step.

\medskip

{\bf 1.2.2. Application to the Mixmaster Universe.} ``Mixmaster Universe''
is defined as the space of solutions  of the
vacuum Einstein equations admitting $SO(3)$ symmetry of the space--like
hypersurfaces (Bianchi IX model, see [Bo]) whose metric
acquires a singularity as near $t\to +0$. The metric
in appropriate coordinates takes the following form:
$$
ds^2=dt^2 -a(t)^2dx^2-b(t)^2dy^2-c(t)^2dz^2. 
$$
\smallskip
The coefficients $a(t),\, b(t),\, c(t)$ are called scale factors.
\smallskip

A family of such metrics satisfying Einstein equations
is given by {\it Kasner solutions}:
$$
a(t)=t^{p_1},\ b(t)=t^{p_2},\ c(t)=t^{p_3};\ \sum p_i=\sum p_i^2=1.
$$
Around 1970, V.~Belinskii, I.~M.~Khalatnikov, E.~M.~Lifshitz
and I.~M.~Lifshitz discovered that most of the
trajectories in Mixmaster Universe exhibit a chaotic behavior
as $t\to +0:$ see [BeKhLi] and subsequent amplifications
in [BoN], [KhLiKSS], [Bar], [May2]. Roughly speaking,
the behaviour of a typical trajectory followed
backwards in time (to the ``Big Bang'')
is described in these papers in the following
way.

\smallskip

Introduce the local logarithmic time $\Omega$  along this trajectory:
$d\Omega:=-\dfrac{dt}{abc}.$ Then $\Omega \to +\infty$
approximately as $-\roman{log}\,t$ as $t\to +0$, and we have:

\smallskip

(i) The time evolution can be divided into
``Kasner eras'' $[\Omega_n,\Omega_{n+1}]$, $n\ge 1.$

\smallskip

(ii) Within each era, the evolution of $a,b,c$ is approximately
described by
Kasner's formula, with variable $p_i$'s which depend
on an additional parameter $u$. If we arrange $p_i$ in the 
increasing order, $p_1<p_2<p_3,$ we have
$$
p_1=-\frac{u}{1+u+u^2},\ p_2=\frac{1+u}{1+u+u^2},\ p_3=\frac{u(1+u)}{1+u+u^2}.
$$
(see [KhLiKSS], formula (2.1).)
The evolution starts with a certain value $u_n>1$, and proceeds as
$u$ diminishes with growing $\Omega$ until $u$ becomes less 
than $1$. After a brief transitional period a new Kasner era starts, with the remarkable
transition formula for the parameter $u_{n+1}$ suggesting
that
continued fractions can be used to model the situation:
$$
u_{n+1}=\frac{1}{u_n-[u_n]}.
$$
\smallskip
(iii) The arrangement of exponents $p_i(u)$ of the scaling
coefficients $a,b,c$ in the increasing order induces generally
a non--identical permutation of these coefficients.
Moreover, during each era several such permutations (Kasner cycles)
occur:
as $u$ diminishes by 1, the old permutation is multiplied
by (12)(3) (see [KhLiKSS], formula (2.3).) When the era
finishes, the permutation (1)(23) occurs
(this is [KhLiKSS], formula (2.2).)

\smallskip

This means that during one era, the largest exponent decreases monotonically,
and governs the same scale factor, $a, b,$ or $c$
which we will call {\it the leading one.} Two other
exponents oscillate between the remaining pair of scaling coefficients.
The number of oscillations is about $k_n:=[u_n].$ Denote
$x_n=u_n-k_n.$

\medskip

Summarizing, we see that in this degree of approximation,
the individual evolution of a typical trajectory is determined
by a number $\alpha\in (0,1)$ whose continued fraction
$[k_1,k_2,k_3,\dots ]$ determines the number of oscillations in each
successive Kasner era. Of course, $\alpha$ is defined only up to a shift,
because the initial point of the backward evolution can
be chosen arbitrarily. Hence the relevant measure 
is the Gauss one.

\smallskip

If we want to keep track of the sequence of the 
leading scale factors as well, we should introduce a set $\bold{P}.$
We claim that in this case it corresponds to the lift of the group
$\Gamma_0(2).$ In fact, consider the action of $GL(2,\bold{Z})$
upon $\bold{P}^1(\bold{F}_2) = \{1,0,\infty\,\}.$ Then
the fractional linear transformation $u\mapsto 1/u$ corresponding
to the transition to the new era, introduces the permutation
(1)(23) of  $\{1,0,\infty\,\},$ whereas the passage to a new cycle 
within one era is described by the transformation
$u\mapsto u-1$ which produces the permutation (12)(3).

\smallskip

Hence the generalized Gauss--Kuzmin theorem leads in this
case to the conclusion that during evolution along
a typical trajectory (i.e. on the set of $\alpha$ of measure 1)
each scale factor becomes the leading one in about
one third of Kasner eras.

\smallskip

(iv) We have not yet connected the proper time $\Omega$
with the variable $u$. This is not directly relevant to our discussion,
but we will do it for completeness, and because
this connection can be beautifully rephrased in terms of ``double--sided
continued fractions'', see [May2], formula (8), [Bar],
formulas (43)--(45), and formula (4.32) below. 

\smallskip

We start with a formula relating the end--points of $\Omega_n$ and $\Omega_{n+1}$ of the $n$--th era with the initial value
$u_n$ inside this era. Namely, introduce one more
parameter $\delta_n>0$ characterizing the relative length
of the era:
$$
\Omega_{n+1}=(1+\delta_n k_n(u_n+1/x_n))\,\Omega_n.
$$
If we put then $\eta_n=(1-\delta_n)/\delta_n$, we have
the following recursion relation:
$$
\eta_{n+1}x_n=\frac{1}{k_n+\eta_nx_{n-1}}.
$$
\smallskip

This means that in terms of the variables $(x_n, y_n:=\eta_{n+1} x_n)$
the transition to the next era is described by invertible
double--sided shift operator
$$ \tilde T: (x,y) \mapsto \left( \frac{1}{x} - \left[ \frac{1}{x}
\right], \frac{1}{y+[1/x]} \right), $$
which is studied in [May2] and [KhLiKSS].

\smallskip
Having  thus completed our
discretized description of the evolution along
an individual trajectory, we have to warn the reader that it refers,
strictly speaking, to {\it another dynamical system} which is defined
on the boundary of a certain compactification
of the phase space of the Mixmaster Universe.
This boundary whose construction involves a nontrivial real blow up
at the $t=0$ subspace was first constructed in [BoN];
see details in [Bo]. The boundary 
is an attractor, it supports an array
of fixed points and separatrices, and the jumps between separatrices
which result from subtle instabilities
account for jumps between 
Kasner's regimes. 
In what sense this picture approximates
the actual trajectories, is not quite trivial question:
cf. the last three paragraphs of the section 2 of
[KhLiKSS], where it is explained that among these
trajectories there can exist
``anomalous'' cases when the description in terms of Kasner
eras does not make sense, but that they are, in a sense, 
infinitely rare.

\medskip

{\bf 1.3. The integral kernel operator.} In this subsection, we
consider the formal operator (1.1) in another functional
space, and show that it admits there a representation with
integral kernel, generalizing that of [Ba], [BaYu].

\smallskip

We will have to assume additionally that the 
$GL(2,\bold{Z})$--set $\bold{P}$ (see 1.1) is such that 
the action
$$ 
t \mapsto \left( \matrix 0&1\\1&k \endmatrix \right) (t)
$$
depends only on $k\, \roman{mod}\, N$ for an appropriate
integer $N$. This assumption is satisfied for instance 
when $G$ is a congruence subgroup. In the following, we fix such $N$.

\smallskip

With this assumption, (1.1) can be written as
$$ 
(L_sf)(z,t)= \sum_{p=-N+1}^0 \, \,  \sum_{k:\,k\ge 1,\,
k\equiv p\, (N)} \frac{1}{(k+z)^{2s}} f\left( \frac{1}{k+z},
\left( \matrix
0 & 1\\
1 & p
\endmatrix \right)\,(t)
\right) \eqno(1.3) 
$$
Our new functional space $\bold{H}$ will consist of functions $f(z,t)$
holomorphic on the
sheets $\{\roman{Re}\,z > -1/2 \} \times t$.

\smallskip

We first recall a useful identity which is the essential ingredient in the
arguments of [Ba], [BaYu], [May1], namely
$$
\sum_{k\geq 1} \frac{1}{(k+z)^{2s}} \exp\left( \frac{-\xi}{k+z}
\right) = \xi^{-s+1/2}\int_0^\infty \eta^{s-1/2} e^{-\eta(z+1/2)} \frac{
J_{2s-1}(2\sqrt{\xi\eta})}{2 \sinh(\eta/2)} d\eta . 
\eqno(1.4)
$$
(This is the formula (111) of [May1] rewritten in a way that
looks more similar to the corresponding formula in [Ba], [BaYu].)

\smallskip

For our purpose, it is useful to consider also the following
corollary of (1.4).

\medskip

\proclaim{\quad 1.3.1. Lemma} 
Let $p\in \{ -N+1,\ldots,0 \}$. We have 
$$ 
\sum_{k:\,k\ge 1,\, k\equiv p\, (N)}
\frac{1}{(k+z)^{2s}}\,\exp \left( -\frac{\xi}{k+z} \right)=
$$
$$
N^{2s-2} \xi^{-s+1/2}\int_0^{\infty}
\eta^{s-1/2}\,\exp\,( {-\eta\, (z+p+N/2)})\,
\frac{J_{2s-1}(2\sqrt {\xi\eta})}{2\,\sinh\,(N\eta /2)}\,d\eta . 
\eqno(1.5)
$$
\endproclaim

{\bf Proof.} We have
$$
\sum_{l\ge 1} \frac{1}{(p+lN+z)^2}\,\exp\,\left(- \frac{\xi}{p+lN+z} \right)=
$$
$$
\frac{1}{N^2}
\sum_{l\ge 1} \frac{1}{(l+(z+p)/N)^2}\,\exp\,\left(- \frac{\xi/N}{l+(z+p)/N}. 
\right).
$$
Plugging (1.4) in and redenoting $N\eta$ as the new $\eta$,
we get (1.5).

\medskip

It is convenient to write the functions $f(z,t)$ in the form
$$ 
f(z,t)= \sum_{j=0}^{|\bold{P}|-1} f_j(z) \delta_j(t), 
$$
with respect to a basis of delta functions.
We have
$$ 
\delta_j\left( \left( \matrix 0&1\\1&p \endmatrix \right)
(t) \right) = \sum_l A_{jl}(p)\, \delta_l(t), $$
where $\bold{A}(p)$ is the matrix representing the permutation of the
sheets.

\smallskip

With this notation, we can write $L_s$ in the form
$$ 
(L_s  f)(z,t)=  \sum_{j=0}^{|\bold{P}|-1} \sum_{p=-N+1}^0
\sum_{k:\,k\ge 
1,\, k\equiv p\, (N)} f_j\left(\frac{1}{k+z} \right) \frac{1}{(k+z)^{2s}}\,
\delta_j(\left( \matrix 0&1\\1&p \endmatrix \right)
(t))=
$$
$$ 
\sum_{j,l=0}^{|\bold{P}|-1} \sum_{p=-N+1}^0 \sum_{k:\,k\ge 1,\,
k\equiv p\, (N)} 
f_j\left(\frac{1}{k+z} \right) \frac{1}{(k+z)^{2s}}
A_{jl}(p)\delta_l(t). 
$$

We now introduce the following operators. Let $\Cal{L}$ denote the
Fourier--Laplace transform
$$ 
(\Cal{L} g)(z,t):=\int_0^\infty e^{-\xi z} g(\xi,t) d\xi. 
$$
We also define the multiplication operator
$$ (T  g)(\xi,t):= e^{-\xi/2} g(\xi,t), $$
as in [Ba], [BaYu], and
$$ 
(S g)(\xi, t):= S(\xi)\,  g(\xi,t),  
$$
with
$$ 
S(\xi)= \frac{(1-e^{-N\xi})^{1/2}}{\xi^{-1/2+s}}.
$$
The operator $\Cal{L}T$ is an isometric isomorphism between the space
$\Cal{L}_2((0,\infty)\times \bold{P})$ and $\bold{H}$ endowed with the norm
$$ 
\| f \|^2 := \frac{1}{2\pi} \sum_{j=0}^{|\bold{P}|-1}
\int_{-\infty}^\infty | f_j(0+iy) |^2 dy. 
$$ 
Define the matrix function
$$
\Theta_{jl}(\xi):= e^{-(N-1)\xi/2}  \sum_{p=-N+1}^0 e^{-\xi p} A_{jl}(p).
\eqno(1.6)
$$
We can and will choose its square root with complex valued
real analytic entries
$\Theta_{jl}^{1/2}(\xi)$:
$$ 
\Theta_{jl}(\xi) = \sum_{k=0}^{| \bold{P} |-1} \Theta_{jk}^{1/2}(\xi)
\Theta_{kl}^{1/2}(\xi),
$$ 
for all $\xi\in (0,\infty)$.
We define
$$
 \hat \Theta_{jl} (\xi):= e^{-(N-1)\xi/4}
\Theta^{1/2}_{jl}(\xi).
\eqno(1.7)
$$
The results of the following Lemma and Proposition represent the analog in
our setting of the main result 
of [Ba], see also 7.4.2 of [May1].

\medskip

\proclaim{\quad 1.3.2. Lemma}
For a function $f(z,t)=\sum_j f_j(z)\, \delta_j(t)$ such that 
$$ 
f_j= \Cal{L} T S^{-1} h_j, 
$$
we can write the operator $L_s$ in the form
$$
(L_s f)(z,t)= \sum_{l=0}^{|\bold{P} |-1} (
\sum_{j=0}^{|\bold{P}|-1} \int_0^\infty e^{-\eta z} S(\eta)^{-1}
\Theta_{jl}(\eta) 
e^{-\eta /2} e^{-N\eta/4}\times
$$
$$
 \int_0^\infty \tilde \kappa(\xi,\eta)
e^{N\xi/4} e^{-\xi/2} h_j(\xi) d\xi 
d\eta  )  \, \, \delta_l(t) 
\eqno(1.8) 
$$
with a function $\tilde \kappa(\xi,\eta)$ satisfying
$\tilde \kappa(\xi,\eta)=\tilde \kappa(\eta,\xi)$. 
\endproclaim
\smallskip

{\bf Proof.} Assume that the functions $f_j(z)$, for $j=0,\ldots, 
|\bold{P} |-1$, are in 
the range of the operator $\Cal{L} T $, namely
$$ 
f_j = \Cal{L} T  g_j, 
$$
for some function $g_j$ in $\Cal{L}_2(0,\infty)$. This means that we can
write 
$$ 
f_j(z)= \int_0^\infty e^{-\xi z} e^{-\xi /2} g_j(\xi) d\xi. 
$$
We now apply (1.5) and obtain
$$ 
= N^{2s-2} \int_0^\infty \xi^{-s+1/2}\int_0^\infty \eta^{s-1/2}
e^{-\eta(z + p + 
N/2)}  \frac{J_{2s-1}(2\sqrt{\xi\eta})}{2 \sinh(N\eta/2)} e^{-\xi/2}
g_j(\xi) d\xi d\eta. 
$$
With our previous definition of $\Theta(\xi)$ as in (1.6)
we can write
$$ 
\sum_p \sum_{k:\,k\ge 1,\, k\equiv p\, (N)}  f_j\left(\frac{1}{k+z} \right)
\frac{1}{(k+z)^{2s}} A_{jl}(p) = 
$$
$$ 
N^{2s-2} \int_0^\infty e^{-\eta z} \int_0^\infty (\frac{\eta}{\xi})^{1/2-s}
\frac{J_{2s-1}(2\sqrt{\xi\eta})}{2 \sinh(N\eta/2)} e^{-(\xi+\eta)/2}
\Theta_{jl}(\eta) g_j(\xi) d\xi d\eta. 
$$ 
Now, for $g_j(\xi)=S(\xi)^{-1} h_j(\xi)$, this can be rewritten as
$$ 
\int_0^\infty e^{-\eta z} S(\eta)^{-1} \Theta_{jl}(\eta)
\int_0^\infty \frac{N^{2s-2} J_{2s-1}(2\sqrt{\xi\eta})}{\sqrt 2
\sinh(N\eta/2)^{1/2}} e^{-(\xi+\eta)/2} e^{-N\eta/4} 
\frac{g_j(\xi)}{\xi^{1/2-s}}  d\xi \, d\eta =
$$
$$ 
\int_0^\infty e^{-\eta z} S(\eta)^{-1} \Theta_{jl}(\eta)
\int_0^\infty \tilde\kappa(\xi,\eta) e^{-(\xi+\eta)/2} e^{-N\eta/4} e^{N\xi /4}
h_j(\xi) d\xi \, d\eta, 
$$ 
where we have set
$$ 
\tilde\kappa(\xi,\eta):= \frac{N^{2s-2} J_{2s-1}(2\sqrt{\xi\eta})}{ 2
\sinh(N\eta/2)^{1/2} \sinh(N\xi/2)^{1/2} }. 
$$

\smallskip

Now the final step.

\medskip

\proclaim{\quad 1.3.3. Proposition}
On the range of 
$\Cal{R}: =\Cal{L}TS^{-1} \hat\Theta$, the operator $L_s$ satisfies
$$ 
\Cal{R}^{-1} L_s \Cal{R} = \Cal{M}, 
$$
where $\Cal{M}$ is the integral kernel operator
$$ 
(\Cal{M} \zeta)(\eta, t)= \sum_{i=0}^{|\bold{P} |-1} \int_0^\infty
\sum_{j=0}^{ | \bold{P} |-1 } \Cal{M}_{ij}(\eta,\xi) \zeta_j(\xi) d\xi \, \,
\delta_i(t),  
$$ 
with $\zeta(\xi,t) = \sum_j \zeta_j(\xi) \delta_j(t)$.
The integral kernel is of the form
$$ 
\Cal{M}_{ij}(\eta,\xi) =  \kappa(\eta,\xi) \sum_r
\hat\Theta_{jr}(\xi) \hat\Theta_{ri}(\eta),  
$$
where the function 
$$ 
\kappa(\xi,\eta) = e^{(N/4 -1/2)(\xi+\eta)} \tilde\kappa(\xi,\eta) 
$$ 
still satisfies $\kappa(\xi,\eta)=\kappa(\eta,\xi)$, but the integral
kernel is in general not symmetric.
\endproclaim

\smallskip

{\bf Proof.} For a function
$$ 
f(z,t)=\sum_{j=0}^{ | \bold{P} | -1} f_j(z) \delta_j(t) =\sum_{j=0}^{|\bold{P} |-1}
\int_0^\infty e^{-\xi z} e^{-\xi/2} S(\xi)^{-1} \sum_{i=0}^{ | \bold{P}
|-1} \hat \Theta_{ij}(\xi) \zeta_i(\xi) d\xi \, \delta_j(t), 
$$
we have
$$ 
(L_s  f)(z,t) = \sum_{l=0}^{|\bold{P} |-1} (L_s  f)_l(z) \delta_l(t), 
$$
with $(L_s  f)_l(z)$ of the form
$$ 
\sum_{j=0}^{ | \bold{P} |-1 } \int_0^\infty e^{-\eta z}
\Theta_{jl}(\eta) e^{-\eta/2} S(\eta)^{-1} 
\int_0^\infty e^{-N\eta/4}  \tilde\kappa(\xi,\eta) e^{N\xi /4} e^{-\xi/2}
\sum_{i=0}^{ | \bold{P} |-1 } \hat\Theta_{ij}(\xi) 
\zeta_i(\xi)  d\xi d\eta. 
$$

We can write this equivalently as
$$ 
\sum_{k,i=0}^{ | \bold{P} |-1 }  \int_0^\infty e^{-\eta z}
\hat\Theta_{kl}(\eta) e^{-\eta/2} S(\eta)^{-1} \int_0^\infty
\Cal{M}_{ik}(\eta,\xi) \zeta_i(\xi)  d\xi d\eta, 
$$ 
where we define
$$ 
\Cal{M}_{ki}(\eta,\xi) := e^{(N/4-1/2)(\eta+\xi)} \tilde\kappa(\eta,\xi) 
\sum_{j=0}^{ | \bold{P} |-1 }\hat\Theta_{ij}(\xi)\hat\Theta_{jk}(\eta). 
$$
Thus, we have obtained
$$ 
L_s \, \Cal{L} T S^{-1} \hat\Theta \, \, \zeta = \Cal{L} T S^{-1} 
\hat\Theta \, \Cal{M}  \, \, \zeta. 
$$
\medskip

{\bf 1.4.  Remarks about the $l$--adic case.} Let $l$ be a prime number,
$\bold{Z}_l$ (resp. $\bold{Q}_l$) the ring of
$l$--adic integers (resp. the field of all $l$--adic numbers).
Put also 
$$
\overline{\Cal{Z}}_l:=\left\{ \frac{n}{l^r}\,|\,n,r\in\bold{Z}, r\geq
0\right\}\, 
\cap\, (0,l) \subset \bold{Q},\ \Cal{Z}_l:= \overline{\Cal{Z}}_l\,\setminus
\,\{1,2,\dots ,l-1\}.
\eqno(1.9)
$$
Every irrational $\alpha\in\bold{Q}_l$ has a unique representation in the form
$$
\alpha = \frac{a_{-r}}{l^r}+\dots +\frac{a_{-1}}{l} +a_0+a_1l+\dots
+a_sl^s+\dots, \  a_i\in\{0,1,\dots ,l-1\}.
$$
The $l$--adic norm of this $\alpha$ is $|\alpha |_l=l^r$ if $a_{-r}\ne 0.$
We then put
$$
[\alpha ]_l:= \frac{a_{-r}}{l^r}+\dots +\frac{a_{-1}}{l} +a_0 \in 
\overline{\Cal{Z}}_l .
\eqno(1.10)
$$
In more invariant terms, $[\alpha ]_l$ is the unique
element in $\overline{\Cal{Z}}_l$ such that $|\alpha - [\alpha ]_l|<1.$

\smallskip

Repeating the usual reasoning, one sees that each irrational
$l$-adic $\alpha$ with $|\alpha |_l<1$ determines a unique sequence
of $k_n(\alpha )\in \Cal{Z}_l$ and $x_n(\alpha )\in
l\bold{Z}_l\setminus \{0\}$ such 
that
$$
\alpha = [\,k_1(\alpha ),\dots ,k_{n-1}(\alpha ), k_n(\alpha )+x_n(\alpha )\,]
\eqno(1.11)
$$
for each $n\ge 1$. 
We get thus the formalism of the theory of continued fractions
in the $l$--adic setting, in which $\bold{R}, \bold{Z}, [0,1)$
are replaced respectively by $\bold{Q}_l,\overline{\Cal{Z}_l}, l\bold{Z}_l$.
Notice that the successive convergents are still rational numbers, but
the incomplete quotients $k_n(\alpha )$ generally are not
integral.

\smallskip

It will be convenient to restrict ourselves to
irrational $\alpha$ in $\bold{Q}_l^*\setminus \bold{Z}_l^*.$
In this case all $k_n(\alpha )$ will belong to $\Cal{Z}_l$.

\smallskip

The shift operator is given by $T:\,\alpha\mapsto \alpha^{-1}-[\alpha^{-1}]_l$,
and it transforms $\bold{Q}_l^*\setminus \bold{Z}_l^*$ into itself.

\smallskip

The definition of a (deformed) transfer operator, however,
presents interesting new problems. We can consider two
basic options.

(A) We can try to define the formal transfer operator by the classical formula
$$
(L_sf)(x):=\sum_{k\in \Cal{Z}_l}
\frac{1}{(x+k)^{2s}}
f \left(\frac{1}{x+k}\right)
\eqno(1.12)
$$
in which $\bold{Z}$ is replaced by $\Cal{Z}_l.$

\smallskip

We could have included a second argument $t$, but did not do it in
order to focus on the peculiarity of (1.12) apparent already in this
straightforward 
version of the classical setting. Namely, we should
not imagine $f$ as a function taking real or complex values. 
In fact, otherwise $\dfrac{1}{(x+k)^{2s}}$
will be defined only at rational points $x$ and will not tend to zero
as $k$ runs over $\Cal{Z}_l$ since from the archimedean viewpoint,
$\Cal{Z}_l$ is a dense set inside $[0,l)$ so that (1.12) will tend
to diverge, unless $f$ is highly discontinuous and tends to zero
when the denominator of the argument tends to infinity.
But this last property will be lost after an application of $L_s$.

\smallskip

However, $\Cal{Z}_l$ is discrete and unbounded in the
$l$--adic sense, so that (1.12) still makes sense as an operator
in various $l$--adic function spaces. For example,
one can consider the space of  analytic functions
on $l\bold{Z}_l$ represented by convergent series $\sum_{n\ge 0}a_nx^n$.
\smallskip

This remark, and parallels with the theories of
$l$--adic uniformization and Drinfeld modules,
suggest the following problems.

\smallskip

(i) Find a natural Banach space of $l$--adic functions
in which (1.12) would define a compact operator.

Since compact operators are nuclear in the $l$--adic
theory, this would allow us to define $l$--adic Selberg's zeta
values at integral points $2s>0$ as $\roman{det}\,(1-L_{2s})$.

\smallskip

(ii) Define the set
of $l$--adic reduced matrices by the same prescription
as (0.5), but this time with $k_i$ running over $\Cal{Z}_l$.

\smallskip

Can one find a characterization of this set similar to
that given in [LewZa1] and reproduced in 0.1? Assuming
we know an $l$--adic zeta,
can one find an Euler product for it similar 
to (3.1) below?

\smallskip

(iii) Again assuming a positive answer to the first question,
is there an eigenfunction with eigenvalue $1$ of $L_1$? 
Can one find its measure--theoretic interpretation
in terms of Mazur's theory of $l$--adic integration?
(See e. g. [Man5], \S 8 and 9). 

\smallskip

At this point, we may notice that the passage 
from $T$ to $L_{1}$ in the classical case implicitly involves  
integration with respect to the additively invariant measure,
since these operators are adjoint via the obvious
bilinear form determined by this integration,
cf. formula (2.19) below.

\smallskip

It is well known that this invariant measure becomes
inadequate for integrating $l$--adic valued functions,
for the simple reason that the smaller is, say, an
$l$--adic ball, the larger is
its invariant measure $l$--adically.

\smallskip

Instead, the $l$--adic integration invented by B.~Mazur for treating
$l$--adic $L$--functions, utilizes finitely
additive functions on open/closed subsets of,
say, $\bold{Z}_l$ which take values in bounded
subsets of finite--dimensional $l$--adic spaces.
Such measures produce linear functionals on the
spaces of functions satisfying the Lipschitz condition.
See [Man5], \S 8 and \S 9, where a more general class
of measures of moderate growth is introduced
and studied as well.

\smallskip

Hence another option, perhaps more natural than the
formal prescription (1.12) is:

\medskip

(B) Define and study the transfer operators on various
spaces of $l$--adic measures, as well as appropriate modifications
of the questions (i)--(iii).

\newpage

\centerline{\bf \S 2. Calculation of certain averages} 

\bigskip

{\bf 2.1. P.~L\'evy's lemma.} Let $f$ be a complex valued function defined on
pairs of coprime integers $(q,q^{\prime})$ such that
$q\ge q^{\prime}\ge 1$ and $f(q,q^{\prime})=O(q^{-\varepsilon})$
for some $\varepsilon >0.$ Put for $\alpha\in (0,1]$
$$
l(f,\alpha )=\sum_{n=1}^{\infty} f(q_{n}(\alpha ),q_{n-1}(\alpha )).
\eqno(2.1)
$$

\smallskip

\proclaim{\quad 2.1.1. Proposition (P.~L\'evy, 1929)} We have
$$
\int_0^1l(f,\alpha )d\alpha = \sum{}^{\prime}\ 
\frac{f(q,q^{\prime})}{q(q+q^{\prime})}.
\eqno(2.2)
$$
Sums and integrals in (2.1), (2.2) converge absolutely
and uniformly.
\endproclaim

\smallskip

A notational convention: prime at the summation sign 
as in (2.2) refers to the
domain $q\ge q^{\prime}\ge 1$, $(q,q^{\prime})= 1$.

\medskip

{\bf Proof.} This Proposition is an immediate consequence
of the following statement.

\smallskip

{\it For any $q\ge q^{\prime}\ge 1$ with $(q,q^{\prime})= 1$
there exists a unique $n\ge 0$ such that one can find
$\alpha\in (0,1]$  with $q_{n}(\alpha )=q,\,q_{n-1}(\alpha )=q^{\prime}.$
Moreover, all such $\alpha$ form a semi-interval of length
$\dfrac{1}{q(q+q^{\prime})}$.}

\smallskip

In fact, assume that  such an $\alpha$ and $n$ exist, 
and let $p_{n}(\alpha )/q_{n}(\alpha ),$ 
$p_{n-1}(\alpha )/q_{n-1}(\alpha )$ be the respective
convergents to $\alpha$, then we have
$$
p_{n-1}(\alpha )\,q_{n}(\alpha )-p_{n}(\alpha )\,q_{n-1}(\alpha )=
(-1)^{n}.
$$ 
Together with the conditions $p_k(\alpha )\le q_k(\alpha )$
this allows us to reconstruct $n$ uniquely by induction and
shows that all $\alpha$ with this property fill the
semi-interval
$$
\frac{p_{n-1}(\alpha )z+p_n(\alpha )}{q_{n-1}(\alpha )z+q_n(\alpha )},\quad
z\in (0,1].
$$
(cf. (0.3)).
Conversely, for any $(q,q^{\prime})$ we can start with
complementing  this line by $p\le q, p^{\prime}\le q^{\prime}$ to a reduced $(2,2)$--matrix with determinant
$\pm 1$, and then produce
the continued fraction for $p/q$ with neighboring convergents
$p/q$, $p^{\prime}/q^{\prime}.$ This proves the lemma.

\medskip

It is often convenient to have the summation domain in (2.2)
extended to all $q\ge q^{\prime}\ge 1.$ One can do this by first
extending $f$ to this domain in the following way:
choose a function $\kappa :\bold{N}\to \bold{Z}$
and a number $t$ and put
$$
F(q,q^{\prime}):= \kappa (d)\,d^{-t}\,f(q,q^{\prime}),
\ d:=\roman{g.~c.~d.}\,(q,q^{\prime}).
\eqno(2.3)
$$
Then
$$
\sum{}^{\prime}\ 
\frac{f(q,q^{\prime})}{q(q+q^{\prime})}=
\zeta (\kappa ,t)^{-1}\,
\sum_{q\ge q^{\prime}\ge 1}\ 
\frac{F(q,q^{\prime})}{q(q+q^{\prime})},
\eqno(2.4)
$$
where $\zeta (\kappa ,t):=\sum_{d\ge 1} \kappa (d)\,d^{-t}.$
(This remark is also contained in [L].)

\smallskip

We will now combine this with the results of [Man1] in order to prove
Theorem 0.2.2.

\medskip

{\bf 2.2. Averaging weighted modular symbols.} In this subsection
we keep the notation explained before the statement
of the Theorem 0.2.2. In particular, modular symbols
refer to the group $\Gamma_0(N).$ We start with the identity (20)
in [Man1]:
$$
\sum_{d/m} \sum_{b=1}^d \int_{\{0,b/d\}} \omega =
(\sigma (m)-c_m)\int_0^{i\infty} \varphi^*(\omega ).
\eqno(2.5)
$$
Here $(m,N)=1,$ $\phi^*(\omega)/dz$ is a cusp form
for $\Gamma_0(N)$ with eigenvalue $c_m$ with respect
to the Hecke operator $T_m$, $\sigma (m)$ is the sum of the divisors
of $m.$

\smallskip

Multiply this identity by $m^{-2-t}$ and sum over all $m$
prime to $N$:
$$
\sum_{m:\,(m,N)=1}\frac{1}{m^{2+t}}
\sum_{d/m} \sum_{b=1}^d \int_{\{0,b/d\}} \omega =
\left[ \sum_{m:\,(m,N)=1} \frac{\sigma (m)}{m^{2+t}}-
L_{\omega}^{(N)}(2+t)\right]\,\int_0^{i\infty} \varphi^*(\omega ).
\eqno(2.6)
$$
Any symbol $\{0,\frac{q^{\prime}}{q}\},\
(q,q^{\prime})=1,$
occurs in the guise $\{0,\frac{b}{d}\}$ for some $d/m$ only when $q$ divides $m$,
and then exactly $\tau (mq^{-1})$ times where $\tau$
is the number of divisors. Hence
the integration path in the left hand of (2.6) can be rewritten
in the following way: 
$$
\sum_{m:\,(m,N)=1}\sum_{q/m} \frac{\tau (mq^{-1})}{m^{2+t}}
\sum_{q^{\prime}\le q,\, (q,q^{\prime})=1} \{0,\frac{q^{\prime}}{q}\}=
$$
$$
\sum_{n:\,(n,N)=1} \frac{\tau (n)}{n^{2+t}}
\sum_{q:\,(q,N)=1}
\frac{\sum_{q^{\prime}\le q,\, (q,q^{\prime})=1} \{0,\frac{q^{\prime}}{q}\}}{q^{2+t}}=
$$
$$
\zeta^{(N)} (2+t)^2\,\left[
\sum_{q:\,(q,N)=1}
\frac{1}{q^{2+t}}\,\sum_{q^{\prime}\le q,\, (q,q^{\prime})=1} \{0,\frac{q^{\prime}}{q}\}\right] .
\eqno(2.7)
$$
Moreover, the first series inside the square brackets in (2.6)
equals 
$$
\zeta^{(N)}(1+t)\,\zeta^{(N)}(2+t).
$$
Hence (2.6) divided by $\zeta^{(N)} (2+t)^2$ can be rewritten as
$$
\sum_{q:\,(q,N)=1}
\frac{1}{q^{2+t}}\,\sum_{q^{\prime}\le q,\, (q,q^{\prime})=1} \int_{\{0,\frac{q^{\prime}}{q}\}}\omega=
$$
$$
\left[ \frac{\zeta^{(N)} (1+t)}{\zeta^{(N)} (2+t)}-
\frac{L_{\omega}^{(N)}(2+t)}{\zeta^{(N)} (2+t)^2}\right]\,\int_0^{i\infty} \varphi^*(\omega ).
\eqno(2.8)
$$

The left hand side of (2.8) can be represented
as the right hand side of (2.2) with the function
$f(q,q^{\prime})$ which vanishes for $q/N$ and otherwise
equals
$$
f(q,q^{\prime})=\frac{q+q^{\prime}}{q^{1+t}}\,\{0,\frac{q^{\prime}}{q}\}.
\eqno(2.9)
$$
Let us define a new function $\tilde{f}(q,q^{\prime})$ by the same
formula (2.9) {\it for all}  relatively prime $(q,q^{\prime}).$
Since $N$ is prime, we have 
$\{0,\frac{q^{\prime}}{q}\}=\{0,i\infty\}$ for $N/q.$
Therefore, writing (2.2) for this $\tilde{f}$, and integrating
$\varphi^*(\omega )$, we get on the left hand side
the same expression as on the left hand side of (0.13).
The right hand side becomes the sum of the right hand side of (2.8)
and 
$$
\sum_{d=1}^{\infty}\frac{\phi (Nd)}{(Nd)^{2+t}}\,\int_0^{i\infty}\varphi^*(\omega)
=\left[\frac{\zeta (1+t)}{\zeta (2+t)}
-\frac{\zeta^{(N)} (1+t)}{\zeta^{(N)} (2+t)}\right]\,
\int_0^{i\infty}\varphi^*(\omega)
\eqno(2.10)
$$
where $\phi$ is the Euler function.
This sum equals the right hand side of (0.16).
This completes the proof.

\medskip

{\bf 2.2.1. Comments and variations.} 
The distribution of
modular symbols was studied by D.~Goldfeld ([Gol1], [Gol2])
who has found
interesting connections between the conjectural asymptotic 
behavior of certain sums involving such symbols and
other number--theoretical problems, e. g. the $abc$--conjecture.
One of Goldfeld's conjectures reads:
$$
\sum \Sb c^2M^2+d^2\le X\\c\equiv 0\,(N)\endSb
\{0,\frac{b}{d}\}_N \sim R(iM)\,X
$$
as $X\to \infty$, where 
the sum is taken over matrices in $\Gamma_0(N)$,
$R(iM):=\int_{iM}^{i\infty}$,
and both sides are considered as functionals on the space of
$\Gamma_0(N)$ cusp forms of weight two.

\smallskip

D.~Goldfeld and C.~O'Sullivan introduced
a class of Eisenstein series twisted by modular symbols and
established their analytic properties. The simplest
series of this kind can be represented as the right hand side of (2.2)
if one chooses for $f$ the following function
(depending on $z,s$ as parameters):
$$
\frac{f(q,q^{\prime})}{q(q+q^{\prime})}:=\chi (q)\, \{0,\frac{q^{\prime}}{q}\}_N\,
\sum_{g \in A_{q,q^{\prime}}}\roman{Im}\,(gz)^s
$$
where $A_{q,q^{\prime}}$ is the set of matrices in $\Gamma_0(N)$
with the second column $(q^{\prime},q)^t$.

\smallskip

Here is another class of quite simple functions $f$ that might produce
interesting specializations of (2.2):
$$
\frac{f(q,q^{\prime})}{q(q+q^{\prime})}:=\frac{\chi_1(q)\,
\chi_2(q^{\prime})}{q^{s_1}
q^{\prime s_2}}.
$$
They lead to some identities involving double
logarithms at roots of unity at the right hand side of (2.2).
As Goncharov has shown in [Gon], relations between 
these numbers can be described in terms of the modular complex for $\Gamma_1(N)$.
This stresses the relevance of the modular symbols in the
study of the distribution of continued fractions.

\smallskip

Our last example is a function that was introduced in [AlZa]:
$$
\frac{f(q,q^{\prime})}{q(q+q^{\prime})}:=
x^{\sum k_j(q/q^{\prime})-1}\,q\,\roman{log}_2\,q .
$$

\medskip

{\bf 2.3. Proof of the Theorem 0.2.1.} We now return to the
notation and conventions explained in the paragraph
around formula (0.13). In particular, $\beta$ is real
irrational.
We start with proving that whenever
the limit (0.13) exists, it does not depend on 
$\alpha$ (independence on $x\in H$ with fixed $\alpha$
is obvious).

\smallskip

We will compare the behavior of (0.13) for two geodesics 
$\Gamma_1,\Gamma_2$ ending at $\beta$.
It suffices to consider the case when $\Gamma_1$
starts from $i\infty$, whereas $\Gamma_2$
starts from some real $\alpha < \beta$.

\smallskip

Denote by $p_n(\beta)/q_n(\beta)=p_n/q_n$ the
convergents to $\beta$. If $n$ is large enough and has
the appropriate parity, the respective convergents will have
the following positions on the real line:
$$
\frac{\alpha +\beta}{2} <\frac{p_{n-1}}{q_{n-1}}<\beta<
\frac{p_n}{q_n}\,.
\eqno(2.11)
$$
Besides, we always have
$$
\left|\frac{p_{n-1}}{q_{n-1}}-\beta\right|>
\left|\frac{p_{n}}{q_{n}}-\beta\right|\,.
\eqno(2.12)
$$
Let $[a,b]$ generally denote the geodesic joining $a$ to $b$.
Define two sequences of points in $H$ by
$$
z_n:=\Gamma_1\cap \left[\frac{p_{n-1}}{q_{n-1}},
\frac{p_{n}}{q_{n}}\right],\
\zeta_n:=\Gamma_2\cap \left[\frac{p_{n-1}}{q_{n-1}},
\frac{p_{n}}{q_{n}}\right] \,.
\eqno(2.13)
$$
From (2.11)--(2.13) it follows that
$$
\frac{1}{2q_nq_{n+1}}<\roman{Im}\,z_n <
\frac{1}{2q_{n-1}q_{n}},
$$
and if moreover $n$ is large enough,
$$
\frac{\theta}{2q_{n+2}q_{n+1}}<\theta\,\roman{Im}\,z_{n+1} <
\roman{Im}\,\zeta_n \le
\frac{1}{2q_{n-1}q_{n}}
$$
where $\theta$ is some fixed constant between $0$ and $1$.
The easiest way to convince oneself of this is to look at
a picture containing all the relevant geodesics.

\medskip

\input epsf
\midinsert
$$\vbox{\epsffile{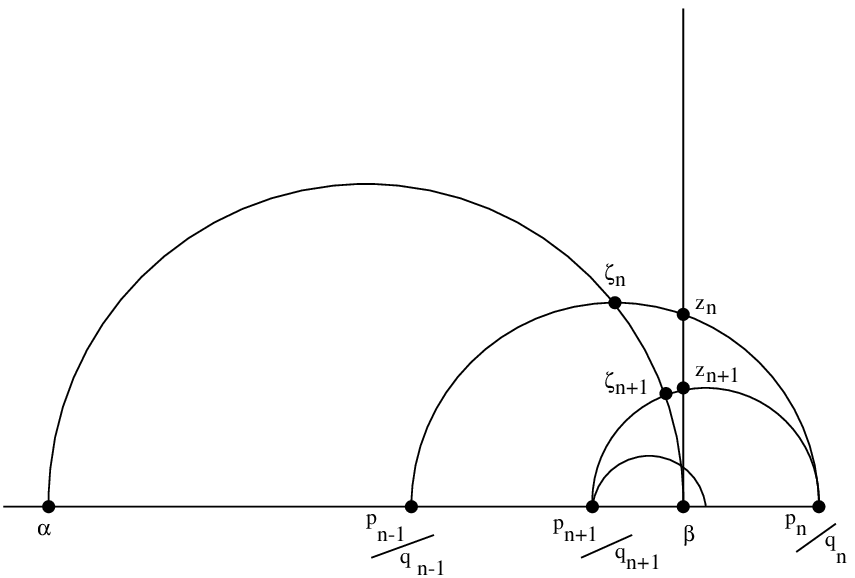}}$$
\endinsert

\medskip

The geodesic distance from any fixed $x_1\in \Gamma_1$
to $z\in \Gamma_1$ equals
$-\roman{log}\,\roman{Im}\,z +O(1)$ as $z\to \beta$. The similar
distance from a fixed $x_2\in \Gamma_2$ to $\zeta \in \Gamma_2$ to $\beta$
equals $-\roman{log}\,\roman{Im}\,\zeta +O(1)$
as $\zeta\to\beta$.

\smallskip

Taking into account our inequalities and the additivity of modular
symbols, we obtain
$$
\frac{1}{T(x_1,z_n)}\,\{x_1,z_n\}=
\frac{1}{T(x_2,\zeta_n)+O(1)}\,[\{x_2,\zeta_n\}+O(1)].
\eqno(2.15)
$$
From this and (2.14) it follows that both limits
exist or otherwise simultaneously, and have a common value whenever 
they both exist.

\smallskip

Moreover, according to the Khintchin--L\'evy theorem
we have for almost all $\beta$
$$
\roman{log}\,q_n = Cn (1+o(1)),\ C=\frac{\pi^2}{12\,\roman{log}\,2}
\eqno(2.16)
$$
as $n\to\infty$.
Hence for almost all $\beta$ we can replace the limit (0.13)
by
$$
\roman{lim}_{n\to\infty}\,\frac{1}{2Cn}\,\{i\infty,z_n\}=
\roman{lim}_{n\to\infty}\,\frac{1}{2Cn}\,\sum_{i=1}^n\left\{\frac{p_{i-1}(\beta)}{q_{i-1}(\beta)},
\frac{p_{i}(\beta)}{q_{i}(\beta)}\right\} .
\eqno(2.17)
$$
Temporarily fixing $n$, we will consider the sum of modular symbols on
the right hand side 
of (2.17) as a function of $\beta$, and then prove that
the resulting sequence of functions weakly converges to zero
in the $\Cal{L}_2$--sense.
For this, we need the following lemma.

\smallskip

\proclaim{\quad 2.3.1. Lemma} Let $\bold{P}_0$ be a finite left
$GL(2,\bold{Z})$--set such that $\roman{Red}^{-1}(t)=\bold{P}_0$
for each $t\in\bold{P}_0.$ Let
$\varphi :\,\bold{P}_0 \to H$
be a function with values in an $\bold{R}$--vector space,
$t_0\in\bold{P}_0.$
Then we have
$$
\roman{lim}_{n\to\infty}\,\frac{1}{n}\,\sum_{i=1}^n \varphi (g_i(x)^{-1}\,t_0)=
\frac{1}{|\bold{P}_0|}\,\sum_{s\in\bold{P}_0} \varphi (s),
\eqno(2.18)
$$
where the limit is taken in the sense of weak convergence in
$\Cal{L}_2([0,1]\times\bold{P}_0)$,
and
$$
g_k(x)=\left(\matrix p_{k-1}(x) & p_{k}(x)\\q_{k-1}(x) & q_{k}(x)
\endmatrix\right) .
$$
\endproclaim
\smallskip

{\bf Proof.} 
Consider the shift 
$$ 
T(x,t)=\left( \frac{1}{x}-\left[ \frac{1}{x} \right], \left(
\matrix -[1/x] &  1\\ 1 & 0 \endmatrix \right)\,(t)
\right) 
$$
We denote by $\lambda$ the measure on $[0,1]\times \bold{P_0}$  given by the standard Lebesgue measure on 
$[0,1]$ and the counting measure on $\bold{P_0}$.

\smallskip

The Gauss--Kuzmin operator $L=L_1$ that we discussed in \S 1 is the
adjoint of this shift $T$, in the  
sense that, for any function $h\in \Cal{L}_1([0,1]\times \bold{P_0},\lambda)$
and any $f\in B_{\bold{C}}$ we have  
$$ 
\int_{[0,1]\times \bold{P_0}} f \cdot Lh \, d\lambda =
\int_{[0,1]\times \bold{P_0}} (f\circ T) \, h \, d\lambda. 
\eqno(2.19)
$$
\smallskip

The eigenfunctional of of $L^*$ denoted $f^*\in K^*$ in \S 1.2
can be taken as
$$
h\mapsto \int_{[0,1]\times \bold{P_0}} h(x,t) d\lambda(x,t).
$$ 
From 1.2 it follows that for any $h\in B_\bold{C}$ we have strong
convergence
$$ 
\lim_{n\to \infty} \frac{1}{n} \sum_{k=1}^{n} (L^k h)(x,t) =
\frac{1}{|\bold{P}_0|\,\roman{log}\,2}\,\frac{1}{1+x}\,
\int_{[0,1]\times \bold{P_0}} h\, d\lambda .
\eqno(2.20) 
$$
According to (2.19), this is equivalent to the convergence
$$
\int_{[0,1]\times \bold{P_0}}\frac{1}{n} \sum_{k=1}^{n}
f(T^k(x,t)) h(x,t) 
d\lambda(x,t)\to
$$
$$
\frac{1}{|\bold{P}_0|\,\roman{log}\,2}\,
\left( \int_{[0,1]\times \bold{P_0}} \frac{f(x,t)}{1+x} 
d\lambda(x,t)\right) \,\int_{[0,1]\times \bold{P_0}} h\, d\lambda  ,
$$
for any $f\in B_{\bold{C}}$ and any test function $h\in B_{\bold{C}}$.
If we consider a function $f(x,t)=\varphi (t)$, independent of $x\in [0,1]$,
we obtain that for any $t$
$$ 
\frac{1}{n} \sum_{k=1}^{n} \varphi (g_k(x)^{-1} t) 
\to \frac{1}{| \bold{P_0} |}\sum_{s\in \bold{P_0}} \varphi(s) 
$$
weakly in $\Cal{L}_2$, because
among the test functions we have all polynomials.

\smallskip

This is equivalent to (2.18) for $\bold{C}$--valued functions
and therefore also for vector--valued ones.

\smallskip

We can now conclude the proof of the Theorem 0.2.1.
If our modular curve is $G_0\setminus\overline{H}$,
we put $\bold{P}_0=PSL(2,\bold{Z})/G_0$ and consider
$\bold{P}_0$ as a left $GL(2,\bold{Z})$--set as explained in
the Introduction.

\smallskip

Since modular symbols are left $G_0$--invariant,
we can find a function $\varphi$ and $t_0\in\bold{P}_0$ such that
$$
\varphi (g_k(\beta )^{-1}\,t_0) =\{g_k(\beta )(0), g_k (\beta )(i\infty )\}=
\left\{\frac{p_{k-1}(\beta)}{q_{k-1}(\beta)},
\frac{p_{k}(\beta)}{q_{k}(\beta)}\right\} .
$$
It follows from (2.18) that the weak limit (2.17) is
$$
\frac{1}{2C\,|\bold{P}_0|}\,\sum_k \{h_k(0),
h_k(i\infty)\}
$$
where $h_k$ now run over a complete set of representatives
of $\bold{P}_0$. But this last sum vanishes. In fact,
let
$$
\sigma =\left( \matrix 0 & -1\\ 1 & 0 \endmatrix \right) .
$$
Then $\{h_k\sigma\}$ as well is a complete system of
representatives, and 
$$
\{\sigma (0),\sigma (i\infty)\}=
-\{0,i\infty\}.
$$
Let us stress that the pointwise behavior of (2.17)
might be wildly oscillating. We proved only that for any measurable set $E$
we have
$$
\roman{lim}_{n\to\infty}\,\frac{1}{n}\,\int_E
\,\sum_{i=1}^n\left\{\frac{p_{i-1}(\beta)}{q_{i-1}(\beta)},
\frac{p_{i}(\beta)}{q_{i}(\beta)}\right\}\,d\beta =0.
\eqno(2.21)
$$

\newpage

\centerline{\bf \S 3. Selberg's zeta function}

\medskip

{\bf 3.1. Notation.} In this section we explain the
definition of Selberg's zeta for subgroups of finite
index $G\subset GL(2,\bold{Z})$ and
$G_0\subset SL(2,\bold{Z})$, their representation
as Fredholm determinant, and relations
to geodesics on modular curves. We closely
follow [LewZa1], pp. 3--6, whose version requires
only minor modifications. For a much more comprehensive
treatment, see [ChMay].

\smallskip

As in [LewZa1], for $g\in GL(2,\bold{Z})$ put 
$$
D(g)=\roman{Tr}\,(g)^2-4\,\roman{det}\,(g),\
N(g)=\left(\frac{\roman{Tr}\,(g)+D(g)^{1/2}}{2}\right)^2,
$$
and call $g$ hyperbolic if  $\roman{Tr}\,(g)$ and $D(g)$
are positive. A hyperbolic matrix
is primitive if it is not a nontrivial power of an element
of $GL(2,\bold{Z})$. For a hyperbolic $g$ set
$$
\chi_s(g) =\frac{N(g)^{-s}}{1-\roman{det}\,(g)\,N(g)^{-1}}.
$$
As above,
put $\bold{P}:=GL(2,\bold{Z})/G$ and denote by
$\rho_{\bold{P}}$ the natural representation of $GL(2,\bold{Z})$
in the space of functions on $\bold{P}$. Finally, put
$$
Z_G(s):=\prod_{ g\in \roman{Prim}}
\prod_{m=0}^{\infty} \roman{det}\,\left[
1-\roman{det}\,(g)^m \,N(g)^{-s-m}\,\rho_{\bold{P}}(g)\right]
\eqno(3.1)
$$
where $\roman{Prim}$ is a set of representatives of all
$GL(2,\bold{Z})$--conjugacy classes of primitive hyperbolic
elements of $GL(2,\bold{Z})$.
 
\smallskip

For $G_0\subset SL(2,\bold{Z}),$ we define $Z_{G_0}(s)$
in the same way, replacing $\roman{Prim}$ by $\roman{Prim_0}$,
a set of representatives of all
$SL(2,\bold{Z})$--conjugacy classes of primitive hyperbolic
elements of $SL(2,\bold{Z})$, and $\bold{P}$ by
$\bold{P}_0$.

\smallskip

\proclaim{\quad 3.2. Theorem} We have
$$
Z_G(s)=\roman{det}\,(1-L_s),\ Z_{G_0}(s)=\roman{det}\,(1-L_s^2)
\eqno(3.2)
$$
where $L_s$ is given by (1.1) and considered as
a nuclear operator in the space $B_{\bold{C}}$.
\endproclaim

\smallskip

We give only a sketch of formal calculations for
$Z_G(s)$. Using notation as in the proof
of Lemma 1.1.2, we have
$$
-\roman{log}\,\roman{det}\, (1-L_s)=
\sum_{l=1}^{\infty}\frac{\roman{Tr}\,L_s^l}{l}=
\roman{Tr}\,\left(\sum_{l=1}^{\infty}\frac{1}{l}
\left(\sum_{n=1}^{\infty}\pi_{s,n}\right)^l\right)=
$$
$$
\roman{Tr}\,\left(\sum_{g\in\roman{Red}} \frac{1}{l(g)}
\pi_s(g)\right) = \sum_{g\in\roman{Hyp}}\frac{1}{k(g)}\,\chi_s(g)\,\tau_g .
\eqno(3.3)
$$
Here we define the operator $\pi_s(g)$ for a reduced matrix $g$ as in (0.5)
as the product of the respective $\pi_{s,k_i}$, and $l(g)$
means its length. $\roman{Hyp}$ denotes a set of representatives
of all conjugacy classes of hyperbolic matrices, and $k(g)$
the maximal integer such that $g=h^{k(g)}.$
The last piece of notation is
$$
\tau_g:=\roman{Tr}\,(\rho_{\bold{P}}(g)) =
\roman{card} \,\{t\in\bold{P}\,|\,g(t)=t\}.
\eqno(3.4)
$$
The appearance of $\tau_g$ is explained by the fact that
our $\pi_s(g)$ acts as the tensor product of the similar operator $\pi_s(g)$
for the case $G=GL(2,\bold{Z})$ and of $\rho_{\bold{P}},$
and our trace is the product of the respective traces.

\smallskip

Using the properties (i)--(iii)  of $\roman{Red}$ summarized 
at the end of 0.1, we can keep rewriting (3.3):
$$
\sum_{g\in\roman{Hyp}}\frac{1}{k(g)}\,\chi_s(g)\,\tau_g=
\sum_{g\in\roman{Prim}}\sum_{k=1}^{\infty}\frac{1}{k}\,
\frac{N(g)^{-ks}\,\tau_{g^k}}{1-\roman{det}\,(g)^k\,N(g)^{-k}}.
\eqno(3.5)
$$
On the other hand,
from (3.1) we find:
$$
-\roman{log}\,Z_G(s)=\sum_{g\in\roman{Prim}}\sum_{m=0}^{\infty}
\roman{Tr}\,
\sum_{k=1}^{\infty} \frac{1}{k}
\,\roman{det}\,(g)^{mk}\,
N(g)^{-(s+m)k}\,\rho_{\bold{P}}(g^k) =
$$
$$
\sum_{g\in\roman{Prim}}
\sum_{k=1}^{\infty} \frac{1}{k}
\frac{N(g)^{-ks}\,\tau_{g^k}}{1-\roman{det}\,(g)^k\,N(g)^{-k}}
\eqno(3.6)
$$
This finishes the formal argument which differs from that
of [LewZa1], \S 1, only by the presence of $\tau_g$. 
The subsequent check of convergence in [LewZa1]
and the argument of \S 3 concerning $SL$ generalize in the same
straightforward way.

\smallskip

Finally, in order to interpret (3.6) in the language of
closed geodesics, it remains only to notice
that if $G\subset GL(2,\bold{Z})$ is the lift of 
$G_0\subset PSL(2,\bold{Z})$ as in the Introduction,
then $G_0\setminus{H}$ can be naturally identified
with $GL(2,\bold{Z})\setminus(H\times\bold{P})$,
and any closed geodesic on the respective
modular curve is covered by the geodesics $[\alpha^-_g,\alpha^+_g]$ 
lying on
those sheets which are left invariant by the respective
hyperbolic matrix $g\in G,$ in agreement with (3.4). 

\newpage

\centerline{\bf \S 4. Non--commutative geometry and the modular
complex}

\medskip

{\bf 4.1. Non--commutative modular curves.} 
As discussed in the Introduction, we want to regard the boundary
$\bold{P}^1(\bold{R})$ with the action of $PSL(2,\bold{Z})$ as a
moduli space of ``non--commutative elliptic curves'', where the
quotient $PSL(2,\bold{Z})\backslash   
\bold{P}^1(\bold{R})$ is itself a non--commutative space. 
According to the general philosophy underlying non--commutative
geometry, this is done by replacing the quotient with the crossed
product of an algebra of functions on $\bold{P}^1(\bold{R})$ by the
action of the group $PSL(2,\bold{Z})$. 
More generally, we can consider the quotients
$G\backslash \bold{P}^1(\bold{R})$
as non--commutative spaces, where $ G$ is a finite index subgroup of
$PSL(2,\bold{Z})$. The classical quotient  
$$ G\backslash (H \cup \bold{P}^1(\bold{Q})), $$
with $H$ the upper half plane, is the modular curve
$G \backslash H$ together with its algebro--geometric
compactification by the set of cusps $G\backslash \bold{P}^1(\bold{Q})$.
The quotient of the full $\bold{P}^1(\bold{R})$ can be regarded as
that part of the analytic boundary which is invisible to the
algebro--geometric compactification, and can be considered as a
``non--commutative modular curve'' when replaced by the crossed
product. We can either consider the crossed product
$$ C(\bold{P}^1(\bold{R}))\rtimes G \eqno(4.1) $$
or, if $\bold{P}$ denotes the coset space $\bold{P}=PSL(2,\bold{Z})/G$, we
can consider the (reduced) crossed product $C^*$--algebra
$$ C(\bold{P}^1(\bold{R})\times \bold{P}) \rtimes PSL(2,\bold{Z}).
\eqno(4.2) $$
The $C^*$--algebras (4.1) and (4.2) are strongly Morita equivalent.

For a discussion of some properties of crossed product
$C^*$--algebras arising from the action of Fuchsian groups on their
limit set, see e.g.~[An--De], [LaSp], [Spi].

We argued in the Introduction that the modular
complex introduced in [Man1] and further studied in [Mer] provides a
definition of cohomology of our boundary space compatible with 
passage to the limit. In this section we show that, in fact, the
modular complex can be related to some standard
homological constructions of non--commutative geometry for the
non--commutative spaces (4.1) or (4.2).

\medskip

{\bf 4.2. Notation.} 
In the following we denote by $\hat X=\bold{P}^1(\bold{R}) \times
\bold{P}$, with $\bold{P}$ the 
coset space $\bold{P}=PSL(2,\bold{Z})/G$.
Moreover, we have
$PSL(2,\bold{Z})=\bold{Z}/2 *\bold{Z}/3$, where we denote by
$$ \sigma:  x\mapsto -1/x \eqno(4.3)  $$
the generator of $\bold{Z}/2$ acting on $\bold{P}^1(\bold{R})$ and by 
$$ \tau: x\mapsto -1/(x-1) \eqno(4.4) $$
the generator of $\bold{Z}/3$.
This action is conjugate to the action on the unit circle by 
rotation by $\pi$ or $2\pi/3$, respectively.
Let $X_{G}=X_{G}(\bold{C})$ denote the modular curve 
$$ X_G=G\backslash \overline{H}, $$
with
$$ \overline{H} =H\cup \bold{P}^1(\bold{Q}). $$

We denote by $\tilde I$ and $\tilde R$ the elliptic points, namely the
orbits $\tilde I= PSL(2,\bold{Z}) \cdot i$ and $\tilde R = 
PSL(2,\bold{Z}) \cdot \rho$. 
We denote by $I$ and $R$ the image in $X_G$ of the elliptic points
$$ I= G\backslash \tilde I \eqno(4.5) $$
$$ R= G\backslash \tilde R, \eqno(4.6) $$
with $\rho=e^{\pi i/3}$.
Finally, for $x$ and $y$ in $\overline{ H }$ we denote by $\langle x,
y \rangle$ the oriented geodesic arc connecting them.

\medskip

\input epsf
\midinsert
$$\vbox{\epsffile{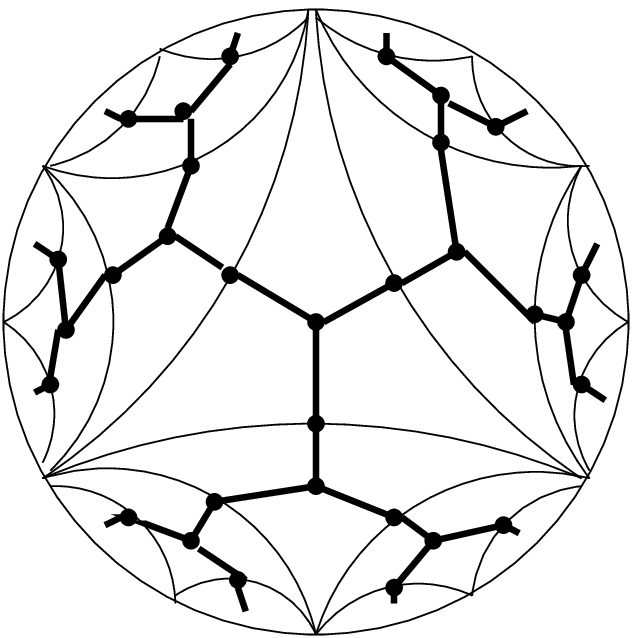}}$$
\endinsert

\medskip

{\bf 4.3. Modular complex}. We consider the following
complex:

0--cells: the cusps $G\backslash \bold{P}^1(\bold{Q})$, and the elliptic
points $I$ and $R$.

1--cells: the oriented half--edges oriented from the parabolic to the
elliptic point:
$$ G \backslash \{ \langle g(i\infty),g i \rangle, g\in
PSL(2,\bold{Z}) \} $$
and the edges
$$ G \backslash \{ \langle g(i), g(\rho) \rangle, g \in
PSL(2,\bold{Z}) \}. $$

2--cells: The images of $E=\{ i,\rho,1+i, i\infty \}$,
$$ G\backslash \{ PSL(2,\bold{Z})\cdot E \}. $$

These cells correspond to the image under the projection 
$\pi: \overline{H} \to X_G$ of all the cells that appear in the figure,
including the vertices on the boundary at infinity of the hyperbolic
disk.

The boundary operators for this complex are given by
$$ \partial : C_2 \to C_1 $$
$$ g E \mapsto g \langle i, \rho \rangle + g \langle \rho, 1+i
\rangle + g \langle 1+i, i\infty \rangle + g \langle i\infty, i
\rangle \eqno(4.7) $$
and
$$ \partial : C_1 \to C_0 $$
$$ g \langle i\infty, i \rangle \mapsto g(i) - g(i\infty) \eqno(4.8) $$
$$ g \langle i, \rho \rangle \mapsto g(\rho)- g(i). \eqno(4.9) $$

The arguments of [Man1] show that this complex computes the homology
of $X_G$,
$$ H_1(X_G)\cong \frac{Ker(\partial: C_1 \to C_0)}{Im(\partial:
C_2\to C_1)}. \eqno(4.10) $$

\medskip

{\bf 4.3.1. Modular complex for relative homology}. There are versions of the 
modular complex considered in [Mer], computing 
relative homology of $X_G$ with respect to cusps and elliptic
points. Here we consider two cases, which differ slightly
from those considered in [Mer]. 
In the modular complex described before, we have
$\bold{Z}[\text{cusps}]=C_0/ \bold{Z}[R\cup I]$. The quotient complex
$$ 0 \to C_2 @>\partial >> C_1 @>\tilde\partial >>
\bold{Z}[\text{cusps}] \to 0, \eqno(4.11) $$ 
where $\tilde\partial$ is the quotient of the boundary operator of the
original modular complex, computes the relative homology $H_1(X_G, 
R\cup I)$. The cycles are given by $\bold{Z}[\bold{P}]$, that is,
combinations of elements $g \langle i, \rho \rangle$, $g$ ranging over
representatives of $\bold{P}$, and by the elements 
$\oplus a_g \langle g(i\infty), g(i) \rangle$
satisfying $\sum a_g g(i\infty)=0$. In fact, these can be represented
as relative cycles in $(X_G,R\cup I)$. 

The subcomplex
$$ 0 \to \bold{Z}[\bold{P}] @>\partial >> \bold{Z}[R\cup I] \to 0,
\eqno(4.12) $$
with $\bold{Z}[\bold{P}]$ generated by the  elements $g \langle i,
\rho \rangle$, computes the homology $H_1(X_G - \text{cusps})$, and
the homology 
$$ H_1(X_G - \text{cusps}, R\cup I) \cong \bold{Z}[\bold{P}]
\eqno(4.13) $$ 
is generated by the relative cycles $g \langle i,
\rho \rangle$.

\smallskip

For convenience of notation, we introduce the same notation used in
[Mer] for the relative homology groups 
$$ H^B_A := H_1(X_G - A, B;
\bold{Z}). $$
These groups are related by the pairing
$$ H^B_A \times H^A_B \to \bold{Z}. \eqno(4.14) $$
In particular, we consider the groups $H^{\text{cusps}}_{R\cup I}$,
$H_{\text{cusps}}^{R\cup I}$, $H^{\text{cusps}}$, and
$H_{\text{cusps}}$.  

We consider an analog of Merel's exact sequences in this setting,
given by the long exact sequence of relative homology 
$$ 0\to H_{\text{cusps}} \to H_{\text{cusps}}^{R\cup I}
@>(\tilde\beta_R,\tilde\beta_I) >> H_0(R)\oplus H_0(I) \to \bold{Z} \to 0,
\eqno(4.15) $$ 
with $H_{\text{cusps}}$ and $H_{\text{cusps}}^{R\cup I}$ as above, and
with 
$$ H_0(I)\cong \bold{Z}[ \bold{P}_I ], \,\,\,\,\,\,
\bold{P}_I= \langle \sigma \rangle\backslash \bold{P} = G\backslash
\tilde I  \eqno(4.16) $$
$$ H_0(R)\cong \bold{Z}[ \bold{P}_R ], \,\,\,\,\,\, \bold{P}_R=
\langle \tau \rangle\backslash \bold{P} = G\backslash \tilde R,
\eqno(4.17) $$ 
that is,
$$ 0\to H_{\text{cusps}} \to \bold{Z}[\bold{P}] \to \bold{Z}[\bold{P}_R]
\oplus \bold{Z}[\bold{P}_I] \to \bold{Z} \to 0. $$ 

\smallskip

We want to compare modular symbols and the non--commutative topology
of the boundary $G\backslash \bold{P}^1(\bold{R})$, in a way that is
compatible with group restrictions $G'\subset G$.

In the case of $H_{\text{cusps}}^{R\cup I}$ and $H_{R\cup
I}^{\text{cusps}}$ the pairing (4.14) gives the
identification of $\bold{Z}[\bold{P}]$ and 
$\bold{Z}^{|\bold{P}|}$, obtained by identifying the elements of
$\bold{P}$ with the corresponding delta functions. 
Thus, we can rewrite the sequence (4.15) as
$$ 0 \to H^{\text{cusps}} \to \bold{Z}^{ | \bold{P} | }
@>(\beta_R,\beta_I) >> \bold{Z}^{ | \bold{P}_I | } \oplus  \bold{Z}^{
| \bold{P}_R |} \to \bold{Z} \to 0, \eqno(4.18) $$
with
$$ H_{R\cup I}^{\text{cusps}} \cong \bold{Z}^{ | \bold{P} | }
\eqno(4.19).  $$    

The map $(\tilde\beta_R,\tilde\beta_I)$ of the relative homology
sequence (4.15) maps $s \mapsto ([s]_R,[s]_I)$,
where $s \in \bold{P}$ corresponds to the generator $g\langle i, \rho
\rangle$, for $g\in PSL(2,\bold{Z})$ the chosen representative of $s\in
\bold{P}$, and $[s]_R \in \bold{P}_R =G\backslash \tilde R$ and $[s]_I
\in \bold{P}_I=G\backslash \tilde I$ are the G--orbits of $g(\rho)$
and $g(i)$, respectively. The map $(\beta_R,\beta_I)$ is given by
$\delta_s \mapsto \delta_{[s]_R}\oplus \delta_{[s]_I}$.

\medskip

{\bf 4.3.2. Algebraic version}. We recall the algebraic formulation of
the modular complex computing $H^{\text{cusps}}= H_1(X_G,
\text{cusps})$, following [Man1] \S 1.8 (a). 

We consider the set of generators $\delta_s$ with $s\in \bold{P}$,
given by the modular symbols $\{ g(0), g(i\infty) \}_G$, with $g$ in
the chosen set of representatives of the cosets $\bold{P}$.
The relations given by the 2--cells can be described as follows.
Consider the subgroup $C$ of $\bold{Z}^{|\bold{P}|}$ with generators
$\delta_s$ and relations $\delta_s \oplus \delta_{\sigma s}$ or 
$\delta_s$ if $s=\sigma s$. Then the homology group $H^{\text{cusps}}$
can be identified with the quotient of $C$ by the subgroup generated
by $\delta_s \oplus \delta_{\tau s} \oplus \delta_{\tau^2 s}$, or
$\delta_s$ if $s=\tau s$. This follows from the arguments of [Man1] \S
1.8 (a). 

\smallskip 

In order to relate this description to the sequence (4.18), consider
first the homology group (4.13), $H_{\text{cusps}}^{R\cup I}=
\bold{Z}[\bold{P}]$. This is generated by the images in $X_G$ of the
geodesic segments $g\gamma_0:=g\langle i, \rho \rangle$, with $g$
ranging over the chosen representatives of the coset space $\bold{P}$. 

Following [Mer], we can identify the dual basis $\delta_s$ of
$H^{\text{cusps}}_{R\cup I}= \bold{Z}^{|\bold{P}|}$ with the images in
$X_G$ of the paths $g\eta_0$, where for a chosen point $z_0$ with
$0< Re(z_0)< 1/2$ and $|z_0|>1$ the path $\eta_0$ is given by the
geodesic arcs connecting $\infty$ to $z_0$, $z_0$ to $\tau z_0$, and
$\tau z_0$ to $0$. These satisfy
$$ [g \gamma_0] \bullet [g \eta_0] =1 $$
$$ [g \gamma_0] \bullet [h \eta_0] =0, $$
for $g G \neq h G$, under the intersection pairing (4.14).

The identification of $H^{\text{cusps}}$, given in terms of generators
and relations as above, with $Ker(\beta_R,\beta_I)$ in the sequence
(4.18) is obtained by the identification 
$$\{ g(0), g(i\infty) \}_G \mapsto g\eta_0,$$ 
so that the relations imposed on the generators
$\delta_s$ by the vanishing under $\beta_I$ correpond
precisely to the relations $\delta_s \oplus \delta_{\sigma s}$ (or 
$\delta_s$ if $s=\sigma s$) and the vanishing under $\beta_R$ gives
the other set of relations $\delta_s \oplus \delta_{\tau s}
\oplus \delta_{\tau^2 s}$ (or $\delta_s$ if $s=\tau s$). 

\smallskip

We shall use this algebraic formulation in the following, when
we relate the group $H^{\text{cusps}}$ to the non--commutative
topology of $G\backslash \bold{P}^1(\bold{R})$.

\medskip

{\bf 4.4. Pimsner exact sequence}. 
We consider the reduced crossed product $C^*$--algebra (4.2). We
recall the setting of Pimsner [Pim] (cf.~[LaSp] for the case of
$C(\bold{P}^1(\bold{R}))\rtimes PSL(2,\bold{Z})$). 
With the notation introduced above, we consider
$$ \Gamma=PSL(2,\bold{Z})=\bold{Z}/2 * \bold{Z}/3$$ 
acting on a tree $T$ with set of
edges $T^1=\Gamma$ and set of vertices $T^0$ given by the cosets
$\Gamma/\Gamma_0$ and $\Gamma/\Gamma_1$, where $\Gamma_0=\bold{Z}/2$
and $\Gamma_1=\bold{Z}/3$.  
This tree $T$ can be realized as a graph in the
2-dimensional hyperbolic space $H$, where the vertices are the
elliptic points $\Gamma \cdot i$ and $\Gamma \cdot \rho$ and the edges
are the geodesic segments $\Gamma \cdot \langle  i,  \rho \rangle$,
represented by bold lines in the figure, giving rise to the
subcomplex (4.12) of the modular complex. 

Associated to this action on a tree, there is a six term exact
sequence [Pim]: 
  
{\eightpoint
$$\diagram
K_0(C(\hat X)) \rto^{\alpha \qquad\qquad}
 & K_0(C(\hat X)\rtimes \Gamma_0)  \oplus K_0(C(\hat X)\rtimes \Gamma_1)
 \rto^{\qquad\qquad \tilde\alpha} &  K_0(C(\hat X)\rtimes \Gamma)  \dto \\ 
K_1(C(\hat X)\rtimes \Gamma) \uto &
K_1(C(\hat X)\rtimes \Gamma_0)  \oplus K_1(C(\hat X)\rtimes \Gamma_1)
\lto^{ \tilde\beta \qquad\qquad} 
 & K_1(C(\hat X))
\lto^{\qquad\qquad \beta} 
\enddiagram  \eqno(4.20) $$ }

We prove the following result that relates the six term exact
sequence (4.20) to the modular complex. 

\smallskip

\proclaim{\quad 4.4.1. Theorem} There exists a natural
isomorphism of the four terms exact sequence 
$$ 0\to Ker(\beta)\hookrightarrow K_1(C(\hat X)) 
@>\beta >> K_1(C(\hat X)\rtimes \Gamma_0)  \oplus K_1(C(\hat X)\rtimes
\Gamma_1) \to Im(\tilde \beta) \to 0  $$
and the exact sequence (4.18),
$$ 0\to H^{\text{cusps}} \to H_{R\cup I}^{\text{cusps}} \to
\bold{Z}^{ | \bold{P}_I | } \oplus  \bold{Z}^{
| \bold{P}_R |} \to \bold{Z} \to 0. $$
These isomorphisms are compatible with the restriction of the group
$G'\subset G$. Moreover, the identification $Ker(\beta)\cong
H^{\text{cusps}}$ is given via the algebraic formulation of \S 4.3.2. 
\endproclaim

\smallskip

{\bf Proof.} First recall that we have natural identifications
$$ K_0(C(\bold{P}^1(\bold{R})) ) \cong \bold{Z} \ \ \ \ \ \ 
 K_1(C(\bold{P}^1(\bold{R})) ) \cong \bold{Z}, $$
given, respectively, by the rank of projections and by the winding
number of the determinant of elements in $GL_n(C(S^1))$. 

Moreover, for the finite groups $\Gamma_j$, there are canonical
isomorphisms,
$$ K^0_{\Gamma_j}(\hat X) \cong
K_0(C(\hat X)\rtimes \Gamma_j) $$
given by
$$ [E] \mapsto [\Gamma(E)], $$
with $E$ a $G$--vector bundles and $\Gamma(E)$ its space of continuous
sections. This gives natural identifications
$$ K_i(C(\hat X)) \cong \bold{Z}^{|\bold{P}|} \ \ \  K_i(C(\hat
X)\rtimes\bold{Z}/2) \cong \bold{Z}^{| \bold{P}_I |} \ \ \ 
K_i(C(\hat X)\rtimes\bold{Z}/3) \cong \bold{Z}^{| \bold{P}_R
|}. \eqno(4.21) $$

Thus, we obtain from (4.20) and (4.21)
$$\diagram 
\bold{Z}^{ |\bold{P}| } \rto^{\alpha
\quad} & \bold{Z}^{ |\bold{P}_I| } \oplus \bold{Z}^{
|\bold{P}_R| } \rto^{\tilde\alpha\quad}  &
K_0(C(\hat X)\rtimes \Gamma)  \dto \\
K_1(C(\hat X)\rtimes \Gamma) \uto &
\bold{Z}^{ |\bold{P}_I| } \oplus \bold{Z}^{
|\bold{P}_R| } \lto^{\quad \tilde\beta} & \bold{Z}^{ |\bold{P}| }
\lto^{\quad \beta} 
\enddiagram \eqno(4.22) $$

The maps in this sequence are defined as in [Pim] \S 1, and they
depend on a choice of fundamental domain for the action of $\Gamma$ on
the tree $T$, which, in our case, is given by the edge
$\langle i, \rho \rangle$ in $T$ and the vertices $\{ i,\rho \}$.

We can split the six term exact sequence (4.22) as
$$ 0\to Ker(\alpha)\hookrightarrow \bold{Z}^{ |\bold{P}| }
@>\alpha>> \bold{Z}^{ |\bold{P}_I|} \oplus
\bold{Z}^{ |\bold{P}_R| } \to Im(\tilde \alpha) \to 0 \eqno(4.23) $$
and
$$ 0\to Ker(\beta)\hookrightarrow \bold{Z}^{ |\bold{P}| }
@>\beta >> \bold{Z}^{ |\bold{P}_I| } \oplus
\bold{Z}^{ | \bold{P}_R | } \to Im(\tilde \beta) \to 0. \eqno(4.24) $$

\smallskip

With the notation of [Pim] \S 4, the morphism $\beta$ (or $\alpha$) in
the Pimsner exact sequence is induced by the maps
$$ \beta_y: C(\hat X) \to C(\hat X)\rtimes \Gamma_{t(y)} \ \ \ \ 
\beta_y (a) = \gamma_{y^t}^{-1} (a) $$
and
$$ \beta_{\bar y} : C(\hat X) \to C(\hat X)\rtimes \Gamma_{o(y)} \ \ \ \  
 \beta_{\bar y}(a) = \gamma_{y^o}^{-1} (a). $$
Here we denote by $y$ the edge $\langle i, \rho \rangle$ in the chosen
fundamental domain for the action of $\Gamma$ on the tree $T$, and
$o(y)=i$, and $t(y)=\rho$ its source and terminus. The groups
$\Gamma_{o(y)}=\bold{Z}/2$ and $\Gamma_{t(y)}=\bold{Z}/3$ are the
stabilizers of these points. Also, here $y^o$
and $y^t$ denote the edges of $T$ with $t(y^t),o(y^o)\in \{ i,\rho
\}$, and $\gamma_{y^t}$ and $\gamma_{y^o}$ are the elements of
$\Gamma$ that satisfy $\gamma_{y^t} y^t =y$ and $\gamma_{y^o} y^o =y$,
as in [Pim] \S 1. The element $\gamma_{y^t}$ acts on
$C(\hat X)$ as the element $\tau$ and $\gamma_{y^o}$ acts as the
element $\sigma$. 

\smallskip

Thus, the morphism 
$$ \beta: K_1(C(\hat X)) \to K_1(C(\hat
X)\rtimes\bold{Z}/2) \oplus K_1(C(\hat X)\rtimes\bold{Z}/3) $$
is precisely the map that sends the generator $\delta_s$, identified
with the homotopy class of determinant functions of elements in
$GL_n(C(\hat X))$ with winding number one around the circle $S^1\times
\{ s \}$ and zero around the circles $S^1 \times \{ t \}$ for $t\neq
s$, to the element $\delta_{[s]_I}\oplus \delta_{[s]_R}$.

We can identify of the maps $\beta$ in (4.24) and
$(\beta_{I},\beta_{R})$  of (4.18). 

This means that there is a natural identification  
$$ Ker(\beta) \cong H^{\text{cusps}} \eqno(4.25) $$
obtained via the algebraic formulation of \S 4.3.2. In fact,
the kernel of $\beta$ in (4.24) can be identified with the subgroup of
$\bold{Z}^{| \bold{P} |}$ of elements $\sum a_s \delta_s$ satisfying 
the relations $a_s + a_{\sigma s} =0$ (or $a_s=0$ if $s=\sigma s$) and
$a_s + a_{\tau s} + a_{\tau^2 s} =0$ (or $a_s=0$ if $s=\tau s$).

The identification (4.25) is compatible with restrictions. In fact, 
suppose given another finite index subgroup $G'$ of $PSL(2,\bold{Z})$,
with $G'\subset G$. This gives a branched cover $X_{G'}@>\pi>> X_G$,
and a surjection $\bold{P}' @>\pi>> \bold{P}$, with $\bold{P}' =
PSL(2,\bold{Z})/ G'$, and a corresponding map
$C(\hat X) \to C(\hat X')$ given by composition with $\pi$. Since the
action of $\tau$ or $\sigma$ on $\bold{P}'$ and $\bold{P}$ commutes
with $\pi$, the maps $\beta_y$ and $\beta_{\bar y}$ are also
compatible with restrictions to $G'\subset G$, and the induced
morphism  $\bold{Z}^{|\bold{P}|} @>\pi^*>> \bold{Z}^{|\bold{P}'|}$
with
$$ \diagram   \bold{Z}^{ |\bold{P}|} \dto^{\pi^*} \rto & \bold{Z}^{
|\bold{P}_I | } \oplus \bold{Z}^{  |\bold{P}_R | } \dto^{\pi^*} \\
\bold{Z}^{ |\bold{P}'|}  \rto & \bold{Z}^{
|\bold{P}_I' | } \oplus \bold{Z}^{  |\bold{P}_R' | }  
\enddiagram $$
is given by $\delta_s \mapsto \oplus_{t\in \pi^{-1}(s)} \delta_t$,
where $\delta_s$ is the homotopy class of determinant functions of
elements in $GL_n(C(\hat X))$ with winding number one around the
circle $S^1\times \{ s \}$, in 
$$ \diagram  K_1(C(\hat X)) \dto \rto^{
\beta\qquad\qquad\qquad} & K_1(C(\hat X)\rtimes \bold{Z}/2)\oplus
K_1(C(\hat X)\rtimes \bold{Z}/3) \dto \\
K_1(C(\hat X')) \rto^{
\beta\qquad\qquad\qquad} & K_1(C(\hat X')\rtimes \bold{Z}/2)\oplus
K_1(C(\hat X')\rtimes \bold{Z}/3) \enddiagram $$

Thus, the identification (4.25) is compatible with restrictions
$G'\subset G$, and the induced map $\pi^* : H^{\text{cusps}} \to
H^{\text{cusps}'}$ in  
{\eightpoint
$$ \diagram 0 \rto &  H^{\text{cusps}} \dto^{\pi^*} \rto & 
\bold{Z}^{|\bold{P}|} \dto^{\pi^*}
\rto^{\beta\quad} &  \bold{Z}^{ |\bold{P}_I
| } \oplus \bold{Z}^{  
|\bold{P}_R | } \rto \dto^{\pi^*} &  \bold{Z} \rto \dto &
0  \\ 
0 \rto &  H^{\text{cusps}'}  \rto & 
\bold{Z}^{|\bold{P}'|} 
\rto^{\beta'\quad} &  \bold{Z}^{ |\bold{P}_I'
| } \oplus \bold{Z}^{  
|\bold{P}_R' | } \rto  &  \bold{Z} \rto  &
0 \enddiagram $$ }
has the following description in terms of modular symbols:
$$ \{ g(0),g(i\infty) \}_G \mapsto \oplus_{t\in \pi^{-1}(s)} \{
g'(0),g'(i\infty) \}_{G'}, $$
with $gG=s  \in \bold{P}$ and $g'G' =t \in \bold{P}'$. This map is the
dual, under the intersection pairing (4.14), to the map $\pi_* :
H_{\text{cusps}'} \to 
H_{\text{cusps}}$ defined by $g\langle i,\rho \rangle \mapsto
g'\langle i,\rho \rangle$ with $g'G' =t \in \bold{P}'$ and $gG= \pi(t)
\in \bold{P}$.

\medskip

The sequence (4.23) differs from (4.24) by the presence of a torsion
term ${\Cal T}$ in $Im(\tilde \alpha)$. The map $\alpha$, which is the
map induced on $K_0$ by the action of $\sigma$ and $\tau$ on $\hat X$,
still satisfies $Ker(\alpha)\cong H^{\text{cusps}}$, but with
different multiplicities, since the morphism for the case $\hat
X=\bold{P}^1(\bold{R})$ is given by
$$ K_0(C(\bold{P}^1(\bold{R}))) @>(2,3)>>
K_0(C(\bold{P}^1(\bold{R}))\rtimes \bold{Z}/2) \oplus
K_0(C(\bold{P}^1(\bold{R}))\rtimes \bold{Z}/3). $$
The torsion term ${\Cal T}$ in $Im(\tilde \alpha)$ depends on the
elliptic elements of $G$. Namely, we have 
$$ Im(\tilde \alpha)= \bold{Z} \oplus {\Cal T} \cong
\bold{Z}^2/\bold{Z}(\ell,1) \oplus T(n_1,\ldots n_k), $$
where the group $G$ has signature $(g;n_1,\ldots,n_k;q)$, with $g$ the
genus and $q$ the number of cusps of $X_G$. Here 
$T(n_1,\ldots n_k)$ is the term computed in [An--De], and $\ell =
l.c.m.(n_1,\ldots n_k)$. 

Thus, from (4.23) and (4.24) we obtain identifications
$$ K_1(C(\hat X)\rtimes \Gamma)\cong H^{\text{cusps}}\oplus
\bold{Z}, \eqno(4.26) $$  
$$ K_0(C(\hat X)\rtimes \Gamma)\cong H^{\text{cusps}}\oplus
\bold{Z}\oplus {\Cal T}, \eqno(4.27) $$
following the identifications of $Ker(\alpha)$, $Im(\tilde\beta)$,
$Ker(\beta)$, and $Im(\tilde\alpha)$ in (4.23) and (4.24).

\bigskip

{\bf 4.5. Cyclic homology}. 
Another way of relating the modular complex and homological
constructions of non--commutative geometry is via cyclic homology. 
We consider the {\it algebraic} crossed product of the algebra
of smooth functions on $\hat X = \bold{P}^1(\bold{R})\times \bold{P}$
by the group $\Gamma=PSL(2,\bold{Z})$,
$$ {\Cal B}:= C^{\infty}(\hat X) \rtimes \Gamma. $$
The $\bold{Z}/2$--graded periodic cyclic cohomology $PHC^*({\Cal B})$
is defines as the 
direct limit over $S: HC^{n-1} ({\Cal B}) \to HC^{n+1} ({\Cal B}) $, where
$S$ is the morphism in Connes' exact sequence relating Hochschild and
cyclic cohomology, [Co] III.1.$\gamma$. It is proved in [Nis] that
there is a six terms exact sequence, analogous to the Pimsner sequence
in $K$--theory, for the  the periodic cyclic cohomology
(or dually for the periodic cyclic homology) of the algebraic
crossed product by a group acting
on a tree. With the notation $A:=C^\infty (\hat X)$ we have 
{\eightpoint
$$\diagram
PHC_0(A) \rto^{\alpha \qquad\qquad}
 & PHC_0(A\rtimes \Gamma_0)  \oplus PHC_0(A\rtimes \Gamma_1) 
 \rto^{\qquad\qquad \tilde\alpha} &  PHC_0(A\rtimes
\Gamma)  \dto \\  
PHC_1(A\rtimes \Gamma) \uto &
PHC_1(A\rtimes \Gamma_0)  \oplus PHC_1(A
\rtimes \Gamma_1) 
\lto^{ \tilde\beta \qquad\qquad} 
 & PHC_1(A)
\lto^{\qquad\qquad \beta} 
\enddiagram  \eqno(4.28) $$ }

Again we can split this six--term exact sequence as {\eightpoint
$$ 0 \to Ker(\beta) \to PHC_1(A) @>\beta>> 
PHC_1(A\rtimes \bold{Z}/2) \oplus PHC_1(A
\rtimes \bold{Z}/3) \to Im(\tilde \beta) \to 0 $$ }
and {\eightpoint
$$ 0 \to Ker(\alpha) \to PHC_0(A) @>\alpha>> 
PHC_0(A\rtimes \bold{Z}/2) \oplus PHC_0(A
\rtimes \bold{Z}/3) \to Im(\tilde \alpha) \to 0. $$ } 
Argument analogous to the case of $K$--theory show that we have an
identification of these sequences with 
$$ 0 \to H^{\text{cusps}} \to k^{|\bold{P}|} \to k^{|\bold{P}_I |}
\oplus k^{|\bold{P}_R|} \to k \to 0,   $$
where we consider homology with coefficients in a field $k= \bold{R}$
or $\bold{C}$, and a corresponding identification
$$ PHC_0({\Cal B})\cong PHC_1({\Cal B}) \cong H^{\text{cusps}} \oplus
k \eqno(4.29) $$  

The relation between the modular complex and
the periodic cyclic homology of ${\Cal B}$ can be derived also via the 
approach of [BN].  

\medskip

{\bf 4.5.1. Groupoids and cyclic homology}. We introduce the groupoid
${\Cal G}_\Gamma$ for 
the action of $\Gamma$ on $\hat X$. This is a Hausdorff locally
compact \'etale groupoid where the morphisms are
$$ {\Cal G}_\Gamma = \hat X\times \Gamma, $$
the objects are ${\Cal G}_\Gamma^0 =\hat X$, 
and the source and target maps $o,t: {\Cal G}_\Gamma \to {\Cal
G}_\Gamma^0$ are given by
$$ o ((x,s),\gamma) = (x,s) \ \ \ \ \ \ t ((x,s),\gamma)= \gamma(x,s),
$$ 
with composition
$$ ((x,s),\gamma)( \gamma(x,s), \gamma')= ((x,s),
\gamma'\gamma).$$
The $C^*$--algebra $C^\infty_c ({\Cal G}_\Gamma)$ of
this groupoid is just the crossed product. 

\smallskip

It is proved in [BN] that the periodic cyclic homology
$PHC_*(C^\infty_c ({\Cal G}_\Gamma))$ is
obtained as a sum of components associated to the torsion conjugacy
classes in $\Gamma$. With the notation of [BN] Corollary 5.9, these
components are given by 
$$ e_{O_\gamma} PHC_n (C^\infty_c ({\Cal G}_\Gamma)) = H_{n + N +
2\bold{Z}} (\hat X^{\gamma} \times_\Gamma E\Gamma),
\eqno(4.30) $$
where $\hat X^\gamma$ is the fixed point set of $\gamma$ acting on
$\hat X$, $N=\dim \hat X^\gamma$, and $\hat X^{\gamma} \times_\Gamma
E\Gamma$ is the homotopy quotient. 

\smallskip

In our case, the only component that contributes in (4.30) 
is the one corresponding to the identity, and, since we consider
homology with coefficients in the field $k$, we can replace the
homotopy quotient by $\Gamma \backslash (\hat X \times H)$, with
$H$ the hyperbolic plane. 
Using the cell decomposition for $\bold{P} \times H$ given by the
modular complex, compatible with the action of $\Gamma$, it is
possible to define
a cell complex computing the homology of the quotient $\Gamma
\backslash (\hat X \times H)$, which recovers the
identification (4.29) via (4.30).

\medskip

{\bf 4.6. Non--commutative geometry and the shift operator}. 
There is another possible way of constructing a ``non--commutative
space'' representing the action of a finite index subgroup $G \subset
PGL(2,\bold{Z})$ on the projective line at infinity
$\bold{P}^1(\bold{R})$, using the shift operator $T$ acting on $[0,1]
\times \bold{P}$, 
$$ T: (x,t) \mapsto \left( \frac{1}{x} - \left[ \frac{1}{x} \right],
\left(\matrix -[1/x] & 1 \\ 1 & 0 \endmatrix\right) (t) \right),
\eqno(4.31), $$ 
as introduced in \S 0.1.1. In fact, the set $X=[0,1]
\times \bold{P}$ meets every orbit of the action of $\Gamma$ on $\hat X =
\bold{P}^1(\bold{R})\times \bold{P}$, and two points $(x,s)$ and
$(y,u)$ in $X$ are equivalent under the action of $\Gamma$ iff there exist
positive integers $m,n$ such that $T^m(x,s)=T^n(y,u)$. 
Thus, it makes sense to consider the action of $T$ on $X$ as a way of
defining a non--commutative analog of the quotient $\hat X/\Gamma$.
The shift $T$ is locally invertible, and it determines  a singly
generated pseudogroup in the sense of 
[Ren], and it defines on $X=[0,1]\times \bold{P}$ an essentially
free singly generated dynamical system, in the sense of Definitions
2.3 and 2.5 of [Ren]. There is an associated  semidirect product
groupoid with arrows
$$ {\Cal G}(X,T)= \{ ((x,t), m-n, (y,s)) | \, T^m (x,t) = T^n (y,s)
\} $$
and objects 
$$ {\Cal G}(X,T)^0 = X \cong \{ ((x,t),0,(x,t) \} \subset {\Cal
G}(X,T).  $$ 
The source and range maps and the multiplication are given
by
$$ p((x,t) ,m-n, (y,s))=(x,t) \ \ \ \ \ \ q((x,t) ,m-n, (y,s))=(y,s)
$$
$$ ((x,t), m-n, (y,s)) \cdot ((x',t'), m'-n', (y',s')) = ((x,t), m+m'
-(n+n'), (y',s')). $$
This is a Hausdorff locally compact \'etale groupoid,
with the topology generated by the basis of open sets
$$ {\Cal G}(T)_{U,V} =\{ ((x,t), m-n, (y,s)) | \, (x,t)\in U, (y,t)\in V, \,
T^m (x,t) = T^n (y,s) \}, $$ 
with $U$ and $V$ open sets where $T^m$ and $T^n$ respectively are
invertible. It is possible to construct a corresponding
$C^*$--algebra $C^*({\Cal G}(X,T))$, which we may also regard as a
non--commutative version of the ``boundary'' $\hat X/\Gamma$. 

\smallskip

It is also interesting to consider the double--sided shift operator
$$ \tilde T : [0,1] \times [0,1] \times \bold{P} \to [0,1] \times
[0,1] \times \bold{P} $$  
$$ \tilde T : (x,y,t) \mapsto \left( \frac{1}{x} - \left[ \frac{1}{x}
\right], \frac{1}{y + \left[ \frac{1}{x} \right]}, \left( \matrix -[1/x] & 1
\\ 1 & 0 \endmatrix \right) (t)  \right), \eqno(4.32). $$
The shift (4.32) is related to the Poincar\'e return map of the
geodesic flow on the modular curve $X_G$, and the one--sided shift
(4.31) is the restriction to the expanding directions (cf.~[ChMay] \S
3.3). The double--sided shift (4.32) is invertible and it defines by
composition an automorphism of the algebra $C(\hat Y)$, with $\hat
Y=[0,1] \times [0,1] \times \bold{P}$. The crossed product 
$C(\hat Y)\rtimes_{\tilde T} \bold{Z}$ gives a natural way of
replacing the set 
of equivalence classes under the action of $T$ by a non--commutative
space. 

\smallskip   

The invariants of the non--commutative geometry of $C(\hat
Y)\rtimes_{\tilde T} \bold{Z}$ should therefore contain some
information on the geodesic flow on the compactified modular curve
$X_G$. For instance the Pimsner--Voiculescu exact
sequence
$$ 0 \to K_1(C(\hat Y)\rtimes_{\tilde T} \bold{Z}) \to K_0(C(\hat
Y)) @>I-\tilde T_*>> K_0(C(\hat Y)) \to K_0(C(\hat Y)\rtimes_{\tilde
T} \bold{Z}) \to 0 $$ 
should be related to the properties of the action of $\roman{Red}$ on
the coset space $\bold{P}$, hence to the dynamical properties of
the geodesic flow.

\medskip

{\bf 4.6.1. Further remarks}.
There are other possible ways of introducing non--commutative geometry at
the boundary of the modular curves. For instance,
if we consider the disconnection of $\bold{P}^1(\bold{R})$ at all the
fixed points of the parabolic elements of $G$, as defined in \S2 of
[Spi], we obtain a totally disconnected compact Hausdorff space
$\Sigma_G$. By the results of [Spi], the crossed product $C^*$--algebra
$C(\Sigma_G)\rtimes G$ is isomorphic to a Cuntz--Krieger algebra
$O_A$, where the matrix $A$ of zeroes and ones corresponds to a
subshift of finite type associated to a choice of the
fundamental domain for the group $G$ as in [BS].
The $K_0$ and $K_1$ of this $C^*$--algebra can be
computed respectively as the cokernel and the kernel of $(I-A^t)$.
These invariants should also contain some information on the
boundary of the modular curves. Using this same technique we
can construct a 
Cuntz--Krieger algebra $O_A$ with the Markov partition determined by
the action of $PSL(2,\bold{Z})$ on $\overline{H} \times \bold{P}$ with
fundamental domain $E \times \bold{P}$, with $E=\{ i,\rho,1+i, i\infty
\}$. By [Spi], this determines a disconnection $\Sigma$
of $\bold{P}^1(\bold{R})$ along $\bold{P}^1(\bold{Q})$, such that the
algebra $C(\Sigma \times \bold{P})\rtimes PSL(2,\bold{Z})$ contain an
image of $O_A$.
Similarly, if we consider the disconnection $\Sigma$ of
$[0,1]$ at all the rational points and the compact totally disconnected
space $X' = \Sigma \times \bold{P}$ with the action of
the shift operator $T$, we obtain a Markov shift as in \S 4 of [Ren],
such that the $C^*$--algebra $C^*({\Cal G}(X',T))$ is a
generalized Cuntz--Krieger algebra for infinite matrices, in the sense
of [EL], with partial isometries $S = \{ (x,1,Tx), x\in U \}$,
with $U$ the sets of the Markov partition.
Again, it should be possible to relate in interesting ways the
calculation of the $K$--theory for this algebra, according to the
techniques of [EL], to the dynamical properties of the shift operator
$T$ and to the boundary of the modular curves.

\bigskip

\newpage

\centerline{\bf References}

\medskip

[AlZa] J.~C.~Alexander, D.~Zagier. {\it The entropy of a certain
infinitely convolved Bernoulli measure.} J. Lond. Math. Soc.,
(2), 44 (1991), 121--134.

\smallskip

[An--De] C.~Anantharaman--Delaroche, {\it $C^*$-alg\'ebres de
Cuntz-Krieger et groupes fuchsiens}. Operator theory, operator
algebras and related topics (Timi\c soara, 1996), 17--35,  
Theta Found., Bucharest, 1997 

\smallskip

[Ba] K.~I.~Babenko, {\it On a problem of Gauss}.
Dokl. Akad. Nauk SSSR, Tom 238 (1978) No. 5, 1021--1024.
(English translation: Soviet
Math. Dokl. Vol.19 (1978) N.1, 136--140.)  

\smallskip

[BaYu] K.~I.~Babenko, S.~P.~Yurev, {\it On a problem of Gauss}.
Sel. Math. Sov. Vol. 2 (1982) N.4, 331--378.

\smallskip

[Bar] J.~D.~Barrow. {\it Chaotic behaviour and the Einstein equations.}
In: Classical General Relativity, eds. W.~Bonnor et al.,
Cambridge Univ. Press, Cambridge, 1984, 25--41.

\smallskip

[BeKhLi] V.~Belinskii, I.~M.~Khalatnikov, E.~M.~Lifshitz.
Adv.~Phys. 19 (1970), 525-551.

\smallskip

[Bo] O.~N.~Bogoyavlenskii. {\it Methods of qualitative theory of dynamical
systems in astrophysics and gas dynamics.} Moscow, Nauka, 1980.

\smallskip

[BoN] O.~N.~Bogoyavlenskii, S.~P.~Novikov. Sov. Phys. JETP,
64:5 (1973), 1475--1491.

\smallskip

[BS] R.~Bowen, C.~Series, {\it Markov maps associated with Fuchsian
groups}, Inst. Hautes Études Sci. Publ. Math. No. 50 (1979), 153--170.

\smallskip

[BN] J.L.~Brylinski, V.~Nistor, {\it Cyclic cohomology of \'etale
groupoids}, $K$--theory 8 (1994) N.4 341--365.

\smallskip

[ChMay] C.-H.~Chang, D.~Mayer. {\it Thermodynamic formalism
and Selberg's zeta function for modular groups}. 
Regular and Chaotic Dynamics, vol. 5, No 3 (2000), 281--312.

\smallskip

[Co] A.~Connes, {\it Noncommutative geometry}, Academic Press, 1994.

\smallskip

[CoDS] A.~Connes, M.~Douglas, A.~Schwarz. {\it Non--commutative
geometry and Matrix theory: compactification on
tori.} Preprint hep--th/9711162.

\smallskip

[EL] R.~Exel, M.~Laca, {\it The $K$--theory of Cuntz--Krieger algebras
for infinite matrices}, K-theory 19 (2000) N.3 251--268.

\smallskip

[Gol1] D.~Goldfeld. {\it Zeta functions formed with modular symbols.}
In: Proc. Symp. Pure Math.. vol. 66, 1 (1999), 111--121.

\smallskip

[Gol2] D.~Goldfeld. {\it The distribution of modular symbols.}
In: Number Theory in Progress (A.~Schinzel's Festschrift,
vol. 2, Walter de Gruyter, Berlin--New York, 1999, 849--865.

\smallskip

[Gon] A.~Goncharov. {\it The double logarithm and Manin's complex 
for modular curves.} Math. Res. Letters, 4(1997), 617--636.

\smallskip

[KhLiKSS] I.~M.~Khalatnikov, E.~M.~Lifshitz, K.~M.~Khanin,
L.~N.~Schur, Ya.~G.~Sinai. {\it On the stochasticity
in relativistic cosmology.} J.~Stat.~Phys., 38:1/2 (1985),
97--114.

\smallskip

[KraLS] M.~Krasnosel'skij, Je.~Lifshits, A.~Sobolev.
{\it Positive linear systems.} Heldermann Verlag, 1989.

\smallskip

[LaSp] M.~Laca, J.~Spielberg, {\it Purely infinite $C^*$--algebras
from boundary actions of discrete groups}, J. reine angew. Math. 480
(1996) 125--139.

\smallskip

[LawZa] R.~Lawrence, D.~Zagier. {\it Modular forms and
quantum invariants of 3--manifolds.} In: Asian J. Math.
3 (1999) (Atiyah's Festschrift), No. 1, 93--107.

\smallskip

[L] P.~L\'evy. {\it Sur les lois de probabilit\'e dont
d\'ependent les quotients complets et incomplets d'une fraction
continue.} Bull.~Soc.~Math.~France, 557 (1929), 178--194.

\smallskip

[LewZa1] J.~Lewis, D.~Zagier. {\it Period functions and
the Selberg zeta function for the modular group.}
In: The Mathematical Beauty of Physics, Adv. Series in
Math. Physics 24, World Scientific, Singapore, 1997, pp. 83--9.

\smallskip

[LewZa2] J.~Lewis, D.~Zagier. {\it Period functions for Maass 
wave forms.} Preprint.

\smallskip

[Man1] Yu.~Manin. {\it Parabolic points and zeta-functions of modular
curves.} Math. USSR Izvestija, vol. 6, No. 1 (1972), 19--64,
and Selected Papers, World Scientific, 1996, 202--247.

\smallskip

[Man2] Yu.~Manin. {\it Explicit formulas for the eigenvalues of Hecke
operators}. Acta Arithmetica, 24 (1973), 239--249.  

\smallskip

[Man3] Yu.~Manin. {\it Three--dimensional hyperbolic geometry as
$\infty$--adic Arakelov geometry.} Inv. Math., 104 (1991), 223--244. 
 
\smallskip

[Man4] Yu.~Manin. {\it Mirror symmetry and quantization of 
abelian varieties.}
Preprint  math.AG/0005143

\smallskip

[Man5] Yu.~Manin. {\it Periods of parabolic forms and $p$--adic Hecke
series.} Math. USSR Sbornik, 21:3 (1973), 371--393,
and Selected Papers, World Scientific, 1996, 268--297.

\smallskip

[Man6] Yu.~Manin. {\it Real Multiplication project and non-commutative
geometry}, Lectures at MPI, May 2001.

\smallskip

[May1] D.~Mayer. {\it Continued fractions and related
transformations.} In: Ergodic Theory, Symbolic Dynamics
and Hyperbolic Spaces, Eds. T.~Bedford et al., Oxford
University Press, Oxford 1991, pp. 175--222.

\smallskip

[May2] D.H.~Mayer, {\it Relaxation properties of the Mixmaster
Universe}, Phys. Lett. A, vol. 121 , nr. 8,9 (1987) 390--394.

\smallskip

[Mer] L.~Merel. {\it Intersections sur les courbes modulaires.}
Manuscripta Math., 80 (1993), 283--289.

\smallskip

[Nis] V.~Nistor, {\it Group cohomology and cyclic cohomology of
crossed products}, Invent. Math. 99 (1990) 411--424.

\smallskip

[O'S] C.~O'Sullivan. {\it Properties of Eisenstein series formed
with modular symbols.} J. reine u. angew. Math.,
518 (2000), 163--186.

\smallskip

[Pim] M.~Pimsner, {\it $KK$--groups of crossed products by groups
acting on trees}, Invent. Math. 86 (1986) 603--634.

\smallskip

[Ren] J.~Renault, {\it Cuntz--like algebras}, Operator theoretical
methods (Timi\c soara, 1998), 371--386,  Theta Found., Bucharest,
2000. 

\smallskip



[Sch] F.~Schweiger. {\it Ergodic theory of fibred systems and metric
number theory.} Clarendon Press, 1995.

\smallskip

[Soi] Y.~Soibelman. {\it Quantum tori, mirror symetry and deformation
theory.} Preprint math.QA/0011162

\smallskip

[Spi] J.S.~Spielberg, {\it Cuntz--Krieger algebras associated with
Fuchsian groups}, Ergod. Th. Dynam. Sys. (1993) 13 581--595.

\smallskip

[Wa] M.S.~Waterman, {\it Remarks on invariant measures for number
theoretic transformations.} Monatshefte f\"ur Mathematik, 79 (1975)
157--163.

\smallskip

[Za] D.~Zagier. {\it Vassiliev invariants and a strange
identity related to the Dedekind eta--function}.
Preprint MPIM 99--78.

\enddocument